\documentclass[reqno,11pt]{amsart}

\usepackage{amsmath, amsthm, a4, latexsym, amssymb}

\def\@rmrk#1#2{\refstepcounter
    {#1}\@ifnextchar[{\@yrmrk{#1}{#2}}{\@xrmrk{#1}{#2}}}

\makeatletter\@addtoreset{equation}{section}\makeatother

 \newfont{\bfit}{cmbxti10 scaled 1200}

 \newcommand{\reals}{{\mathbb{R}}}
 \newcommand{\eps}{\varepsilon}

 \newcommand{\R}{\mathbb{R}}
 \newcommand{\N}{\mathbb{N}}
 
 \newcommand{\prob}{\mathbb{P}}

 \newcommand{\me}{\mathbb{E}}
 \newcommand{\Q}{\mathbb{Q}}
 
 \newcommand{\1}{{\sf 1}}

 \newcommand{\skric}{{\mathcal C}}

 \newcommand{\skrii}{{\mathcal I}}

 \newcommand{\skrim}{{\mathcal M}}

 \newcommand{\skris}{{\mathcal S}}
 \newcommand{\skrit}{{\mathcal T}}

 \newcommand{\skrix}{{\mathcal X}}
 \newcommand{\heap}[2]{\genfrac{}{}{0pt}{}{#1}{#2}}

\renewcommand{\subsection}{\secdef \subsct\sbsect}
\newcommand{\subsct}[2][default]{\refstepcounter{subsection}
\vspace{0.15cm}
{\flushleft\bf \arabic{section}.\arabic{subsection}~\bf #1  }
\nopagebreak\nopagebreak}
\newcommand{\sbsect}[1]{\vspace{0.1cm}\noindent
{\bf #1}\vspace{0.1cm}}

{\nopagebreak {\hfill\rule{2mm}{2mm}}\\ }

\newtheorem{theorem}{Theorem}[section]
\newtheorem{lemma}[theorem]{Lemma}
\newtheorem{cor}[theorem]{Corollary}

\newtheoremstyle{thm}{1.5ex}{1.5ex}{\itshape\rmfamily}{}
{\bfseries\rmfamily}{}{2ex}{}

\newtheoremstyle{rem}{1.3ex}{1.3ex}{\rmfamily}{}
{\itshape\rmfamily}{}{1.5ex}{}
\theoremstyle{rem}

\refstepcounter{subsubsection}

\def\thebibliography#1{\section*{Bibliography}
  \list%
  {\arabic{enumi}.}
    {\settowidth\labelwidth{[#1]}\leftmargin\labelwidth
    \advance\leftmargin\labelsep
    \parsep0pt\itemsep0pt
    \usecounter{enumi}}
    \def\newblock{\hskip .11em plus .33em minus .07em}
    \sloppy                   
    \sfcode`\.=1000\relax}

\begin{document}
{\Large \bf Large deviations of Markov chains indexed by random trees}\\

\hspace{0.5cm}AMIR DEMBO,\footnote{Research
partially supported by NSF grant \#DMS-0072331.}
{\it Stanford University}\\
\hspace*{0.5cm}PETER M\"ORTERS, {\it University of Bath}\\
\hspace*{0.5cm}SCOTT SHEFFIELD, {\it Microsoft Research}\\

\thispagestyle{empty}
\vspace{0.2cm}

\vspace{0.5cm}

\centerline{\bf \small Abstract} \begin{quote}{\small Given a finite
typed rooted tree $T$ with $n$ vertices, the
{\em empirical subtree measure}
is the uniform measure on the $n$
typed subtrees of $T$ formed by taking all descendants of a single vertex.
We prove a large deviation principle in $n$, with explicit rate function,
for the empirical subtree measures of multitype Galton-Watson trees conditioned
to have exactly $n$ vertices. In the process, we extend
the notions of shift-invariance and specific relative entropy---as
typically understood for Markov fields on deterministic graphs such as
$\mathbb Z^d$---to Markov fields on random trees.  We also
develop single-generation empirical measure large deviation principles
for a more general class of random trees
including trees sampled uniformly
from the set of all trees
with $n$ vertices. }
\end{quote} \vspace{0.5cm}

\begin{tabular}{lp{13cm}}
{\it Keywords:} & Tree-indexed Markov chain, branching Markov chain,
random tree, Galton-Watson tree, multitype Galton-Watson process,
multitype Galton-Watson tree, marked tree, large deviation principle,
empirical pair measure, empirical offspring measure, process level.\\
\multicolumn{2}{l}
\end{tabular}

\section{Introduction}\label{intro}

The {\em empirical measures} of Markov fields on large, deterministic subsets $\Lambda$ of $\mathbb
Z^d$---and the limit points of these empirical measures---play a central role in statistical
physics and the theory of Gibbs measures.  The limit points are always shift-invariant, and the
rate functions of the empirical measure large deviation principles are generally defined in terms
of {\em specific relative entropy} or {\em specific free energy},
see, e.g., Chapters 14--16 of \cite{G88}.

When $\mathbb Z^d$ is replaced with a random graph, the large deviation analysis of even the
simplest models---say, Ising or Potts models---becomes more difficult.  How does one even define
``shift-invariance,'' for example, when the graphs on which the models are defined are random and
almost surely possess no translational symmetries?  What is the most natural analog of ``specific
relative entropy''?  For that matter, what is the most useful definition of ``empirical measure''?

The purpose of this paper is to answer the above questions for some natural
random planar rooted tree models.
By {\em planar} we mean that the offspring of each vertex are implicitly
ordered---from left to right; this ordering determines an embedding of the tree in the plane.

Given a finite
planar rooted tree $T$ with $n$ vertices
with
types drawn from a finite type set $\skrix$, the {\em empirical subtree measure}
$\nu^T$ is the uniform measure on the $n$ typed subtrees of $T$ that are formed
by taking all descendants of a single vertex of
$T$.  We will prove a large deviation principle, with an explicit rate function
defined in terms of specific relative entropy on the empirical subtree measures
of multitype Galton-Watson trees conditioned to have exactly $n$ vertices.

The rate function of this large deviation principle will be infinite on measures that lack a
natural ``shift-invariance'' property.  A shift-invariant measure $\nu$ on trees may be either
almost surely finite or almost surely infinite. In either case, we will show that every
shift-invariant measure can be ``extended backwards'' to describe the ``infinite past'' of a sample
from the tree. We may also view this {\em backward tree} construction as a general technique for
examining the steady state of a randomly expanding system.  It is on these backward tree measures
that we will actually define specific relative entropy, as the conditional entropy of the offspring
measure at the root {\em given} its infinite past.

One motivation for pursuing this problem is the study of {\em tree-indexed
Markov chains}, defined as follows.  First we sample a tree from some probability
measure,
and then, given this tree, we run a Markov chain on the vertices of the tree in such a way that the
state of a vertex depends only on the state of its parent. The result of this two-step experiment
can also be interpreted as a \emph{typed tree}. We always look at probabilities with respect to the
whole experiment, or, in the language of random environments, at the \emph{annealed} probabilities.
These tree-indexed process are a natural concept of increasing interest in probability and
applications (see, e.g., \cite{BP94}, \cite{Pe95} and \cite{LPP95}), often as a new way of looking
at existing models.  Our analysis will show that large deviations results, which are well-known for
classical Markov chains, can be extended to Markov chains indexed by random trees.

When we restrict our attention to a single generation of the empirical measure (the ``empirical
offspring measure'') or to a type of empirical measure on typed edges (the ``empirical pair
measure'') we will obtain a generalized large deviation principle for which the classical Markov
results (as developed in, e.g., \cite{DZ98} and the references therein) are a special case.  In
fact, these turn out to be among the rare problems for which large deviation rates can be stated
completely explicitly in a closed form.  Indeed, the rates we find in this setting are hardly more
complicated than the rates for classical Markov chains.
For example, our rate functions are simple enough to allow one to
compute the pressure and related macroscopic quantities for
Gibbs measures corresponding to a short-range potential with
configuration space that is the set of all
typed rooted trees of $n$ vertices with types in $\skrix$.
This  is in sharp contrast with the large deviation principle for
the distance from the root of simple random walk on
supercritical Galton-Watson trees, for which no explicit rate function
is known,
see \cite{DGPZ02}.

In another application, from the case of binary trees and
uniform distribution of types, we calculate an
explicit growth rate for the total number of binary trees of size $n$ (odd) with types in a finite
alphabet $\skrix$, which have an empirical pair measure in a given set of measures. In \cite{KM02}
the analogous combinatorial formula for the number of tuples of length $n$ with a given empirical
pair measure was used to analyse the tail behaviour of Brownian intersection local times. We hope
that the formulas derived here give rise to a similar analysis of the tail behaviour of integrated
super-Brownian excursion, as formulas for high moments of intersection local times involve
summation over large binary trees, see e.g.~\cite{LG99}.

%

There are a number of technical issues that make the analysis of tree-indexed Markov chains more
complicated than the analogous work for classical Markov chains.  One arises from the fact that,
for some models of Galton-Watson trees, the probability of having
exactly $n$ vertices is zero for~$n$ in an infinite subset of $\mathbb Z$.
It is therefore necessary to restrict our attention to
those ~$n$ for which the probability is
positive and to prove lower bounds on probabilities that apply
only for select values of~$n$.  Another arises from the possibility of
an unbounded number of offspring at a single step, which necessitates the
use of a technical ``mass exchange'' argument in Lemma~\ref{apx}.

The precise statements of our results are given in Section \ref{statement}
beginning with empirical pair and empirical offspring measures and then progressing
to the empirical subtree measures.
The former results will apply to a larger class of random trees than the
latter, which will only be proved for bounded-offspring
multitype
Galton-Watson trees.  The proofs of all of
these results are then given in Section \ref{LDPproof}.

\section{Statement of the results} \label{statement}

By $\skrit$ we denote the set of all finite rooted planar trees $T$,
by $V=V(T)$ the set of all vertices and by $E=E(T)$ the set of
all edges oriented away from the root, which is always denoted
by $\rho$. We write $|T|$ for the number of vertices in the tree $T$,
with the $k$-th generation of $T$ being the subset of vertices of $T$
of distance $k$ from its root and the height of $T$ is the largest
$k$ such that the $k$-th generation of $T$ is non-empty.

Suppose that $T$ is any finite tree and we are given an initial
probability measure $\mu$ and a Markovian transition kernel
$Q:\skrix\times \skrix\to[0,1]$ on a finite alphabet $\skrix$.
We can obtain a \emph{tree indexed Markov chain} $X:V \to \skrix$ by
choosing $X(\rho)$ according to $\mu$ and choosing $X(v)$, for each
vertex $v\not=\rho$, using the transition kernel given the value
of its parent, independently of everything else. If the tree is
chosen randomly, we always consider $X=\{X(v)\, : \, v\in T\}$
under the \emph{joint law} of tree and chain. It is sometimes
convenient to interpret $X$ as a \emph{typed tree}, considering
$X(v)$ as the \emph{type} of the vertex $v$.

We first look at the class of Galton-Watson trees, where the
number of children $N(v)$ of each $v \in T$ is an independent
random variable, with the same law 
$p(\,\cdot\,)=\prob\{N(v)=\,\cdot\,\}$ for all
$v \in T$, such that $0<p(0)<1$. With each finite tree and sample path $X$
we associate a probability measure on
$\skrix\times\skrix$, the \emph{empirical pair measure}~$L_X$, by
$$L_X(a,b)= \frac{1}{|E|} \sum_{e\in E}
\delta_{(X(e_1), X(e_2))}(a,b),
\, \mbox{ for } a,b\in\skrix,$$
where $e_1, e_2$ are the beginning and end vertex of the edge $e\in E$
(so $e_1$ is closer to $\rho$ than $e_2$).
Our first result is a large deviation principle for $L_X$, conditional
upon the event $\{|T|=n\}$ with $n$ chosen such that the latter has
positive probability. 
For its formulation recall the definition of the relative entropy
$H(\cdot\, \| \, \cdot)$ from \cite[(2.1.5)]{DZ98} and Cram\'er's rate function
\begin{equation}\label{Ip-def}
I_p (x) = \sup_{\lambda \in \reals} \, \Big\{ \lambda x -
 \log \Big[ \sum_{n=0}^\infty  p(n) e^{\lambda n} \Big] \Big\} \,,
\end{equation}
as in \cite[(2.1.26)]{DZ98}.
\begin{theorem}\label{main}
Suppose that $T$ is a Galton-Watson tree,
with offspring law $p(\cdot)$ such that $0<p(0)<1-p(1)$,
$\sum_\ell \ell p(\ell) = 1$ and $\ell^{-1} \log p(\ell) \to -\infty$.
Let $X$ be a Markov chain indexed by $T$ with arbitrary initial distribution
and an irreducible Markovian transition kernel $Q$. Then, for $n\to\infty$,
the empirical pair measure $L_X$, conditioned on $\{|T|=n\}$
satisfies a large deviation principle
in the space of probability vectors on $\skrix \times \skrix$
with speed $n$ and the convex, good rate function
\begin{align}
I(\mu)=\left\{ \begin{array}{ll} H(\mu \, \| \,\mu_1\otimes Q)
+ \displaystyle\sum_{a \in \skrix} \mu_2(a)
\, I_p \Big(\frac{\mu_1(a)}{\mu_2(a)}\Big)
 & \mbox{ if $\mu_1\ll \mu_2$,} \\
\infty & \mbox{ otherwise,}
\end{array} \right.
\label{Idef}
\end{align}
where $\mu_1$ and $\mu_2$ are the first and second marginal of $\mu$ and
$\mu_1\otimes Q (a,b)=Q\{ b \, | \, a \} \mu_1(a)$.
\end{theorem}
{\bf Remarks:}

\noindent{$\bullet$} Throughout the paper we implicitly assume that the
conditioning events $\{|T|=n\}$ are of positive probability, that is,
our large deviation approximation of probabilities hold for those values of
$n$ where $\prob\{|T|=n\}>0$. 
For the general structure of the set $S$ of
admissible values, see the proof of Lemma~\ref{subexp}.

\noindent{$\bullet$} In case $\sum_\ell \ell p(\ell) \neq 1$ note that
the distribution of $T$ conditioned on $\{ |T| = n \}$ is exactly
the same as when the offspring law is $p_\theta(\ell) =
p(\ell) e^{\theta \ell}/\sum_j p(j) e^{\theta j}$,
regardless of the value of $\theta \in \reals$. With $0<p(0)<1-p(1)$
there exists a unique $\theta_*$ such that
$\sum_\ell \ell p_{\theta_*} (\ell) =1$.
Hence, Theorem \ref{main} still applies, 
using $I_{p_{\theta_*}}$ in place of $I_p$ in \eqref{Idef}.

\noindent{$\bullet$} The representation \eqref{Idef} of $I(\cdot)$ provides
the interpretation of the large deviations of $L_X$ as the result of two
independent contributions: when $\mu_1=\mu_2$ we have only
the term $H(\mu \, \| \,\mu_1 \otimes Q)$ which
is the rate function for the large deviation principle of
empirical pair measures of the Markov chain with kernel $Q$,
see e.g.~\cite[Section 3.1.3]{DZ98}, while the hard constraint of $\mu_1=\mu_2$
of the Markov chain setting is replaced here by the additional term
$\sum_{a} \mu_2(a) I_p (\mu_1(a)/\mu_2(a))$ reflecting the large
deviations contribution due to the geometry of the tree $T$.

{\bf Examples:}

The class of Galton-Watson trees conditioned on the total size
appears in the combinatorial literature, see e.g.~\cite{MM78}, under the name
\emph{simply generated trees} and  is surveyed in \cite{Al91}. We look at some
interesting examples.

\noindent{$\bullet$} Choose the offspring law $p(\,\cdot\,)$ such that
$p(k)=1-p(0)=1/k.$ In this case $\prob\{|T|=n\}>0$ if and only if $n-1$ is
divisible by $k$. The law of $T$ conditional on $\{|T|=n\}$ is exactly
the same as sampling the tree uniformly from the
collection of all possible $k$-ary trees with $n$ vertices. We have that
$I_p(x)= (x/k) \log x + (1-x/k) \log((1-x/k)/(1-1/k))$, leading to the good
rate function
\begin{align*}
I(\mu)=\left\{ \begin{array}{ll} H(\mu \, \| \,\mu_1\otimes Q)
+ \frac{k-1}{k} H\big( \frac{1}{k-1}(k\mu_2-\mu_1) \, \| \, \mu_2\big)
+ \frac{1}{k} H(\mu_1 \,\| \,\mu_2) & \mbox{ if $k\mu_2\ge \mu_1$,} \\
\infty & \mbox{ otherwise,}
\end{array} \right.
\end{align*}
for the large deviation principle of $L_X$.


\noindent{$\bullet$} Choose the offspring law $p(\,\cdot\,)$ as
the standard Poisson distribution, $p(\ell)=e^{-\ell}/\ell!$ for
$\ell=0,1,2,\ldots$.
Now $\prob\{|T|=n\}>0$ for all $n \geq 1$ and the law of $T$ conditioned
on $\{|T|=n\}$ is that of a tree chosen uniformly from all unordered
trees with $n$ vertices. We have $I_p(x)= 1-x+x \log x$, and
get a large deviations rate of
$I(\mu)=H(\mu \, \| \,\mu_1\otimes Q) + H(\mu_1 \, \| \mu_2)$
in (\ref{Idef}).

\noindent{$\bullet$} Choose the offspring law $p(\,\cdot\,)$ as
$p(0)=p(1)=\cdots=p(k)=1/(k+1)$. Note that this law is only critical
if $k=2$, and recall the second remark following Theorem~\ref{main}.
Again $\prob\{|T|=n\}>0$ for all $n \geq 1$, and now the law of $T$
conditional on $\{|T|=n\}$ is the same as sampling the tree uniformly
from the collection of all ordered trees with $n$ vertices and offspring
number bounded by $k$. {\nopagebreak\hfill$\diamond$\\ }

The result extends to other classes of trees, indeed one can go much
beyond the present setting and consider trees and types chosen simultaneously
according to a \emph{multitype Galton-Watson tree}.
In this situation, in order to obtain more explicit rate functions,
it is useful to replace the empirical pair measure
by a more inclusive object, the \emph{empirical offspring measure}.

We write $\skrix^*=\bigcup_{n=0}^\infty \{n\} \times \skrix^n$ and
equip it with the discrete topology. Note that the offspring of any
vertex $v\in T$ is characterized by an element of $\skrix^*$ and that
there is an element $(0,\emptyset)$ in $\skrix^*$ symbolizing lack of
offspring. For each typed tree $X$ and each vertex $v$ we denote by
$$C(v)=(N(v),X_1(v),\ldots,X_{N(v)}(v))\in \skrix^*$$ the
number and types of the children of $v$, ordered
from left to right.
To each sample chain $X$ we associate a probability measure $M_X$ on
$\skrix\times\skrix^*$ called the \emph{empirical offspring measure},
which is defined by
$$M_X(a,c)=\frac 1{|T|} \sum_{v\in V} \delta_{(X(v),C(v))} (a,c).$$

We now describe the joint law of a tree $T$ and tree-indexed chain $X$,
which defines a multitype Galton-Watson tree. The ingredients
are a probability measure $\mu$ on $\skrix$, serving as the initial
distribution, and an offspring transition kernel $\Q$ from $\skrix$
to $\skrix^*$. We define the law $\prob$ of a tree-indexed
process $X$ by the following rules:
\begin{itemize}
\item The root $\rho$ carries a random type $X(\rho)$
chosen according to the probability measure $\mu$ on $\skrix$.
\item For each vertex with type $a\in \skrix$ the offspring number
and types are given independently of everything else, by the
offspring law $\Q\{ \, \cdot\, \,|\, a\}$ on $\skrix^*$. We write
$$\Q\{ \, \cdot\, \,|\, a\}
=\Q\{ (N,X_1,\ldots, X_N)\in \cdot \,\, | \, a\},$$
i.e. we have a random number $N$ of offspring particles with types
$X_1,\dots,X_N$.
\end{itemize}
We assume that the exponential moments $\Q\{ \, e^{\eta N} \,|\, a\}<\infty$,
for all $a\in\skrix$ and $\eta>0$.
We also need a weak form of irreducibility assumption.
Denote, for every $c=(n,a_1,\ldots,a_n)\in \skrix^*$
and $a\in\skrix$, the \emph{multiplicity} of the symbol~$a$ in $c$ by
$$m(a,c)=\sum_{i=1}^n \1_{\{a_i=a\}}.$$
Define the matrix $A$ with
index set $\skrix\times\skrix$ and nonnegative entries by
$$A(a,b)=\sum_{c\in \skrix^*} \Q\{ c \, | \, b\} m(a,c),
\mbox{ for } a,b\in\skrix,$$
i.e.~$A(a,b)$ are the expected number of offspring of type~$a$
of a vertex of type~$b$.
With $A^*(a,b) = \sum_{k=1}^\infty A^k(a,b) \in [0,\infty]$ we say that
the matrix $A$ is \emph{weakly irreducible} if $\skrix$ can be partitioned
into a non empty set $\skrix_r$ of \emph{recurrent states} and a
disjoint set $\skrix_t$ of \emph{transient states} such that
\begin{itemize}
\item $A^*(a,b)>0$ whenever $b \in \skrix_r$, while
\item $A^*(a,b)=0$ whenever $b \in \skrix_t$ and either $a=b$ or $a \in \skrix_r$.
\end{itemize}
For example, any \emph{irreducible} matrix $A$ has
$A^*$ strictly positive, hence is also weakly irreducible with $\skrix_r=\skrix$.
The multitype Galton-Watson tree is called weakly irreducible (or irreducible)
if the matrix $A$ is weakly irreducible (or irreducible, respectively)
and the number $\sum_{a \in \skrix_t} m(a,c)$ of transient offspring is uniformly
bounded under $\Q$.

Note that a weakly irreducible matrix has $A(a,b)=0$ whenever $b \in \skrix_t$
and $a \in \skrix_r$. Moreover $\skrix_t$ may be ordered such that $A(a,b)=0$
when $a \geq b$ are both in $\skrix_t$. Consequently, the non-zero eigenvalues of a
weakly irreducible matrix $A$ are exactly those of the irreducible
matrix obtained by its restriction to $\skrix_r$.
Recall that, by the Perron-Frobenius theorem,
see e.g. \cite[Theorem 3.1.1]{DZ98}, the largest eigenvalue of an irreducible
matrix is real and positive. Obviously, the same applies to
weakly irreducible matrices.  The multitype Galton-Watson tree
is called \emph{critical} if this eigenvalue is $1$ for the matrix $A$.

Our second main result is a large deviation principle for $M_X$ if
$X$ is a multitype Galton-Watson tree. For its formulation denote,
for every probability measure $\nu$ on $\skrix\times\skrix^*$, by $\nu_1$
the $\skrix$-marginal of $\nu$. We call $\nu$ \emph{shift-invariant}
if $$\nu_1(a)= \sum_{(b,c)\in \skrix\times\skrix^*} m(a,c) \nu(b,c)
\mbox{ for all }a\in\skrix.$$
We denote by $\skrim(\skrix\times\skrix^*)$ the space of probability measures
$\nu$ on $\skrix\times\skrix^*$ with $\int n \, \nu(da \, ,
dc)<\infty$, using the convention $c=(n,a_1,\ldots,a_n)$.
We endow this space with the smallest topology which makes the
functionals $\nu\mapsto \int f(b,c)\, \nu(db\, ,dc)$ continuous,
for $f:\skrix\times\skrix^*\to\R$ either bounded, or
$f(b,c)=m(a,c)\1_{b_0}(b)$ for some $a,b_0 \in\skrix$.
Define the function $J$ on
$\skrim(\skrix\times\skrix^*)$ by
$$J(\nu)=\left\{ \begin{array}{ll} H(\nu \,\|\,  \nu_1\otimes\Q) &
\mbox{ if $\nu$ is shift-invariant,} \\ \infty & \mbox{ otherwise.}
\end{array} \right.$$
In general, the topology on $\skrim(\skrix\times\skrix^*)$ is stronger
than the weak topology, making the function $J$ lower semicontinuous, as
shown in Lemma \ref{rate-lsc}.

\begin{theorem}\label{general}
Suppose that $X$ is a
weakly
irreducible, critical multitype Galton-Watson
tree with an offspring law whose exponential moments are all finite,
conditioned to have exactly $n$ vertices. Then, for $n\to\infty$,
the empirical offspring measure~$M_X$ satisfies a large deviation
principle
in
$\skrim(\skrix\times\skrix^*)$
with speed $n$ and
the convex,
good rate function $J$.
\end{theorem}

\noindent
{\bf Examples:}

{$\bullet$}
The 
situation of Theorem~\ref{main} corresponds to offspring kernels
$\Q\{ \, \cdot \, | a\}$ choosing offspring numbers according to the
law $p(\,\cdot\,)$ and then choosing the offspring types independently,
according to the marginal law $Q\{ \, \cdot \, | \, a\}$ on $\skrix$.
Consequently, Theorem~\ref{main} follows
by contraction from Theorem~\ref{general}, see Section~\ref{s-main} for
more details. As its proof reveals, Theorem \ref{main} applies even when
the law of offspring numbers $p(\,\cdot\,|\,a)$ depends on the type of the parent,
provided the matrix $Q\{ b \,|\, a\} \sum_\ell \ell p(\ell \,|\,a)$ is
weakly irreducible, with largest eigenvalue one (then, of course,
$I_{p(\cdot\,|\,a)}$ replaces $I_p$ in \eqref{Idef}).

{$\bullet$} For a more concrete example contained in our framework, we suppose that
individuals in a population may have two genetic types, $A$ and $B$.
Individual of type $A$ (resp.~$B$) breed offspring according to the law
$p_A$ (resp. $p_B$), typically of the same type, but independently,
mutations occur  with a small probability $p>0$. Denote by $\eta$
the ratio of the mean offspring number of $p_A$ and $p_B$, representing the
genetic advantage of type $A$. In a large family of size $n$ the probability
that the ratio of the numbers of individuals of type $A$ and $B$ in the population
is close to $x\in[0,1]$ is approximately equal to $\exp(-nI(x))$ for
$$I(x)=\inf \Big\{ \frac{x}{x+1} H\big(\nu_A \, \big\| \, q_A \big) +
\frac{1}{x+1} H\big(\nu_B \, \big\| \, q_B\big) \Big\},$$
where $q_{A}(n,m)=p_{A}(n+m)\,\left( \heap{n+m}{m} \right)\,p^m(1-p)^n$
and $q_B(n,m)=p_{B}(n+m)\,\left( \heap{n+m}{m} \right)\,p^n(1-p)^m$
and the infimum is over all probability measures $\nu_A, \nu_B$  on
$\N\times\N$ satisfying
$$x=\sum_{n,m=0}^\infty nx\nu_A(n,m)
+ n \nu_B(n,m) \mbox{ and }
1=\sum_{n,m=0}^\infty mx\nu_A(n,m)
+ m\nu_B(n,m).$$
This rate function is zero exactly at the typical
ratio, which is given by the solution $x>0$ of the equation
$x/(1+x) = (x \eta (1-p)+p)/(x \eta+1)$.
Our result gives the probability of a significant deviation from this ratio,
the precise rate is depending of course
on the exact offspring laws of particles of either genetic type,
represented by $p_A, p_B$. {\nopagebreak\hfill$\diamond$\\ }

We conclude with the extension to a \emph{process level} large deviation
principle. For the rest of this section we assume that the offspring numbers
generated by the kernel $\Q$ are uniformly bounded by some $N_0\in\N$.
We denote by $\bar{\skrix}$
the set of all finite or infinite rooted, planar trees such that every vertex
has at most $N_0$ offspring, with types from the finite alphabet $\skrix$
attached to the vertices. Recall that the fact that the trees are embedded
in the plane imposes an ordering
(say 
from left to right) on the children of each vertex.

The laws of multitype Galton-Watson trees are probability measures on
$\bar{\skrix}$. We equip $\bar{\skrix}$ with the topology
generated by the functions $f:\bar{\skrix}\to\R$ depending only
on a finite number of generations.

If $v\in V$ is a vertex of a tree $T$ and $X\in\bar{\skrix}$ a sample
chain on this tree, we denote by $X^v$ the sample chain obtained from the
subtree of $T$ consisting of
$v$ and all successors of $v$. To each \emph{finite}
sample chain $X$ we associate a probability measure $T_X$ on $\bar{\skrix}$,
the \emph{empirical subtree measure}, which is defined by
$$T_X(x)=\frac 1{|T|} \sum_{v\in V} \delta_{X^v} (x),
\mbox{ for }x\in\bar{\skrix}.$$
To formulate a large deviation principle for the random variable $T_X$ we need
further notation. We denote by
$N[k]$ the number of vertices in generation $k$, and in particular by
$N=N[1]$ the number of children of the root in $T$.
Suppose that $\mu$ is a probability measure on $\bar{\skrix}$ with
$\int N d\mu=1$. Then we can define a \emph{shifted} probability measure
$S(\mu)$ on $\bar{\skrix}$ by
\begin{equation}
\label{s-def}
S(\mu)(\Gamma)= \int d\mu(X) \sum_{i=1}^N  \1_{\{X^{v_i} \in \Gamma\}},
\mbox{ for any Borel set } \Gamma \subset \bar{\skrix}\;,
\end{equation}
where $v_1,\ldots,v_N$ are the children of the root. We call $\mu$
\emph{shift-invariant} if $S(\mu)=\mu$.

To any shift-invariant measure $\mu$ on $\bar{\skrix}$ we can associate
a \emph{backward tree measure} $\mu^*$ in the following way.
Suppose that $X$ is a sample chain on a (finite or infinite) tree of height at
least $k$, and mark a vertex in generation $k$
of $X$ as the \emph{centre} of the tree.
Denote by $\skrix[k]$ the set of all objects
$(x,\zeta)$ (typed tree $x$ with centre at $\zeta$) arising in
this way, endowed with the canonical topology inherited from $\bar{\skrix}$.
For $k\ge l$ there are canonical projections $p_{kl}:\skrix[k] \to \skrix[l]$
obtained by keeping the same centre and removing all vertices from the tree
whose last common ancestor with the centre lived before generation $k-l$.
Note that the root of the
projected tree $p_{kl}X$ is the ancestor of the centre in
generation $k-l$. The spaces $\skrix[l]$ and
projections $p_{kl}, k\ge l$ form a
projective system. Hence there exists a projective limit space
$\underline{\skrix}$, the space of \emph{backward trees},
and canonical projections $p_k:\underline{\skrix}\to\skrix[k]$.
See \cite[Appendix~B]{DZ98} for more information about projective limits.

If $\mu$ is a shift-invariant measure then we can associate a measure $\mu_k$
on $\skrix[k]$ by
$$\mu_k(\Gamma)=\int d\mu(X) \, \sum_{i=1}^{N[k]} \1_{\{(X,v_i) \in \Gamma \}}, \,
\mbox{ for any  Borel set } \Gamma \subset \skrix[k],$$
where $v_1,\ldots,v_{N[k]}$ are the vertices in generation $k$ of $X$.

Shift-invariance of $\mu$ ensures that all $\mu_k$ are
probability measures and that $\mu_l= \mu_k \circ p_{kl}^{-1}$ for all $k\ge l$.
Hence, by Kolmogorov's extension theorem, there exists a unique probability
measure $\mu^*$ on $\underline{\skrix}$ such that $\mu^*\circ p_k^{-1}=\mu_k$.
This is the backward tree measure $\mu^*$ associated to $\mu$.

For each $k \geq 1$
we denote by $\mathfrak{p}_{1,k}:\skrix [k] \to\skrix [k]$
the projection obtained by removing all vertices of distance at least $k+2$
from the root and all those of distance $k+1$ from the root whose parent is
to the right of the centre.  Similarly, we denote by
$\mathfrak{p}_{0,k}:\skrix [k] \to\skrix [k]$ the projection which in addition
to all the vertices removed by $\mathfrak{p}_{1,k}$ also removes all
children of the centre. Note that
$p_{kl}\circ \mathfrak{p}_{0,k}=\mathfrak{p}_{0,l} \circ p_{kl}$
and $p_{kl} \circ \mathfrak{p}_{1,k}=\mathfrak{p}_{1,l} \circ p_{kl}$
for all $k \geq l$. Hence, the projective limits
$\mathfrak{p}_1:\underline{\skrix}\to\underline{\skrix}$
and $\mathfrak{p}_{0}:\underline{\skrix}\to\underline{\skrix}$ of
$\mathfrak{p}_{1,k}$ and $\mathfrak{p}_{0,k}$, respectively, are
well defined with $p_k\circ \mathfrak{p}_0=\mathfrak{p}_{0,k} \circ p_k$
and $p_k \circ \mathfrak{p}_1=\mathfrak{p}_{1,k} \circ p_k$
for all $k \geq 1$ (heuristically, $\mathfrak{p}_1$
is the projection obtained by removing all vertices of the backward tree
further from the root than the centre except the children of the
centre and those of the vertices to the right of the centre whose
distance from the root is the same as the centre, with
$\mathfrak{p}_0$ removing also the children of the centre).
If $\Q$ is an offspring transition kernel, we
define $\mu^*\circ \mathfrak{p}_0^{-1}\otimes \Q$ as the probability measure
generated by starting with a backward tree sampled according to
$\mu^*\circ \mathfrak{p}_0^{-1}$ and adding independently
offspring according to $\Q$ to the centre.
Let $\skrim(\bar{\skrix})$
be the set of probability measures on $\bar{\skrix}$.
Define the function $K$ on $\skrim(\bar{\skrix})$ by
$$K(\mu)=\left\{ \begin{array}{ll} H\big(\mu^*\circ
\mathfrak{p}_1^{-1} \, \big\| \,
\mu^*\circ \mathfrak{p}_0^{-1}\otimes \Q\big) &
\mbox{ if $\mu$ is shift-invariant,} \\ \infty & \mbox{ otherwise.}
\end{array} \right.$$
We equip $\skrim(\bar{\skrix})$ with the smallest topology which
makes the functionals $\mu\mapsto \int f\, d\mu$ continuous,
for each continuous and bounded $f:\bar{\skrix}\to\R$.

\begin{theorem}\label{process}
Suppose that $X$ is an irreducible, critical multitype Galton-Watson
tree with
uniformly bounded offspring sizes,
conditioned to have exactly $n$ vertices. Then, for $n\to\infty$,
the empirical subtree measure~$T_X$ satisfies a large deviation
principle
in $\skrim(\bar{\skrix})$
with speed $n$ and
the convex,
good rate function $K$.
\end{theorem}


We now give a brief overview over the following sections, which contain the
proofs of our results. First we need to establish the fact that for a critical
multitype Galton-Watson tree our conditioning events $\{|T|=n\}$ decay with an
exponential rate zero over the set of admissible values of $n$. The proof of
this fact, well-known for single-type Galton-Watson trees, requires a careful
analysis of the lattice structure of the set $S=\{ n\in\N \, : \,
\prob\{|T|=n\}>0\}$ in the multitype case, and is of some independent interest.
This result is proved in Section~\ref{s-decay}.

Equipped with this result, in Section~\ref{s-upper}
the upper bound of Theorem~\ref{general} is derived. Exponential
tightness is established in the topology on $\skrim(\skrix\times\skrix^*)$ using
the moment conditions imposed on $\Q$. Based on the exponential Chebyshev inequality
we first represent the upper bound in a variational form, and then solve the
variational problem. Nonstandard arguments arise in the proof from the fact that
we endow $\skrim(\skrix\times\skrix^*)$ with a topology, which is stronger than
the weak topology of measures. This is necessary in order to make the set of
shift-invariant measures a closed set in our topology.

The lower bound, proved in Section~\ref{s-lower}, is based on a change of measure
technique. As we allow for potentially unbounded offspring numbers intricate
approximation arguments are needed to show that this change of measure provides
sufficient freedom to represent a sufficiently large class of offspring measures.
In Section~\ref{s-main} we prove Theorem~\ref{main} by contraction
from Theorem~\ref{general}.

Finally, in Section~\ref{s-process} we prove Theorem~\ref{process}. For this
purpose we first extend Theorem~\ref{general} from one-generation offspring
measures to $k$-generation offspring measures, see Lemma~\ref{kstep}.
This extension is based on expanding the statespace and needs crucially the fact
that in Theorem~\ref{general} we are only requiring \emph{weak} irreducibility.
The step from $k$-generation offspring measures to empirical subtree measures is
then based on the Dawson-G\"artner Theorem.

\section{Proof of the large deviation principles} \label{LDPproof}

\subsection{On the rate of decay of $\prob\{|T|=n\}$.}\label{s-decay}

An important role in our proofs is played by the fact that
for critical multitype Galton-Watson trees the probability
$\prob\{|T|=n\}$ decays only subexponentially on the set $S$ of
integers $n$ where the probability is positive. We exclude
the trivial case when $S$ fails to be infinite from our consideration
(in particular, we assume throughout that $\mu(\skrix_r)>0$).

\begin{lemma}\label{subexp}
Suppose $T$ is the random tree generated by a
weakly
irreducible,
critical multitype Galton-Watson~tree
with finite second moment. Then
$$\lim_{\heap{n\to\infty}{n \in S}}
\frac 1n \log \prob\{ |T|=n\}=0.$$
\end{lemma}

\begin{proof}
Recall that the number of children of any given $v \in T$ with
types in $\skrix_t$ is uniformly bounded. Moreover, if
$X(u) \in \skrix_t$ for some $u \in T$ then there are only
types from $\skrix_t$ in the
sample chain $X^u$ consisting of $u$ and all successors of $u$,
and the height of the corresponding subtree $T^u$ is uniformly
bounded (by the size of $\skrix_t$). Let
$G(v) = \sum_i |T^{u_i}|$ over the
children $u_1,u_2,\ldots$ of $v$ such that $X(u_i) \in
\skrix_t$.
Hence $G(v)$ is also uniformly bounded, say by
$N_1<\infty$. For $c \in \skrix^*$ let $c|\skrix_r$ be the natural restriction
of $c$ to $\skrix_r^*$. For each $b \in \skrix_r$, $c \in \skrix_r^*$
and $g \in \{0,\ldots,N_1\}$
let $\Q_r \{ (c,g) \,|\, b \}$ denote the probability
induced by $\Q$ that given $X(v)=b$ we have
$C(v)|\skrix_r^*=c$ and $G(v) = g$. Then, for each $c_r \in \skrix_r^*$,
$$\sum_{g=0}^{N_1} \Q_r \{ (c_r,g) \,|\, b \} = \sum_{\{c \in \skrix^*
\,: \, c|{\skrix_r} = c_r \}} \, \Q \{ c \,|\, b \}\,,$$
so $\Q_r$ is a transition probability measure from $\skrix_r$ to
$\skrix_r^* \times \{ 0, \ldots, N_1\}$ such that
$A_r(a,b)=\sum_{c,g} m(a,c) \Q_r \{ (c,g) \,|\, b \}$ is exactly
the restriction of the matrix $A$ to $\skrix_r$. In particular, since
$A$ is weakly irreducible and critical, it follows that $A_r$ is
irreducible and critical on $\skrix_r$. Further, $\Q_r$ constructs
the restriction of the multitype Galton-Watson tree $X$ to $\skrix_r$
with $G(v)$ keeping track of the number of vertices
with types in $\skrix_t$ that have been omitted as a result of being
in $T^u$ for some child $u$ of $v$ such that $X(u) \in \skrix_t$.
Thus, fix a type $a\in\skrix_r$
and construct a multitype Galton-Watson tree
with law $\prob$, for $\mu=\delta_a$ as follows: Start at size $n=0$
with one \emph{active} vertex $\rho$ of type $a$. At each future step
choose an active vertex $v$ uniformly from all active vertices,
independently of everything else, provide it with offspring $C(v)$
according to $\Q_r \{ \cdot \, | \, X(v)\}$, adding $G(v)+1$
to the current tree size $n$, deactivating $v$ and activating its offspring.
When there are no active vertices left, the process
terminates, producing the restriction to $\skrix_r$ of
a typed tree of law $\prob$ and size $n$ for $\mu=\delta_a$.

Let $p_{a,b}(n)$ be the probability that when the size is $n$ we have exactly
one active vertex, 
which is of type $b$. For any $a_1, a_2, a_3\in\skrix_r$ and
positive integers $n_1, n_2$ we have
\begin{equation}\label{onestar}
p_{a_1,a_2}(n_1)p_{a_2,a_3}(n_2)\le p_{a_1,a_3}(n_1+n_2).
\end{equation}
Indeed, $p_{a_1,a_2}(n_1)p_{a_2,a_3}(n_2)$ is the probability of
having exactly one active vertex when the size is $n_1$ and again when the size
is $n_1+n_2$,
having types $a_2$ and $a_3$, respectively.

Since the restricted multitype Galton-Watson tree is
irreducible, starting with $a \in \skrix_r$ active
vertices of each type appear with positive probability and our
procedure allows each active vertex to eventually remain the only
active vertex with positive probability. Hence for any
$a_1,a_2 \in \skrix_r$,
there exists $n$ such that $p_{a_1,a_2}(n)>0$.
Together with \eqref{onestar} this suffices to make the structure
of the sets $$S_{a,b}=\big\{n\in\N \, : \, p_{a,b}(n)>0\big\}$$
for $a,b \in \skrix_r$,
analogous to that of the sets $\{n\in\N \, : \, (P^n)_{a,b}>0\big\}$
for a finite state irreducible Markov chain with transition matrix~$P$.
Namely, there exists a period $d=$gcd~$S_{a,a}$, independent of
$a \in \skrix_r$,
and $k_{a,b}\in\{0,\ldots,d-1\}$ such that $S_{a,b}\subset k_{a,b}+d\N$
with $|(k_{a,b}+d\N)\setminus S_{a,b}|<\infty$, see for example the proof
in~\cite[Lemmas 5.5.3, 5.5.4 and 5.5.6]{Du96}. Analogously to the theory of
$d$-periodic finite state irreducible Markov chains, \eqref{onestar}
and subadditivity imply the existence of $I<\infty$ such that,
for all
$a,b \in \skrix_r$,
$$\lim_{l\to\infty} -\frac1{ld}\log p_{a,b}(k_{a,b}+ld)=I.$$
(Indeed, one can take first $a=b \in \skrix_r$ showing
existence of limits $I_{a,a}<\infty$, then show that $I_{a,a}\le
I_{b,b}$ for all $a,b\in\skrix_r$, hence for each
such $a$ and $b$ the limit
$I_{a,b}$ exists and is equal to $I_{a,a}$ by a sandwich argument).
Now let $p_a(n)=\prob\big\{ |T|=n \, | \, X(\rho)=a \big\}$,
$S_a=\{n \, : \, p_a(n)>0\}$ and $\skrix_g=\big\{b \, : \,
\Q_r \{((0,\emptyset),g) \,|\, b\}>0\big\}$, noting that
the latter set is nonempty for some $g$
(otherwise no finite trees are possible).
The event $\{|T|=n\}$ corresponds to \emph{one} active vertex from
$\skrix_g$ at size $n-1-g$ producing
$g$ omitted vertices of types from $\skrix_t$ and
no offspring with type in $\skrix_r$.
Summing over the possible types of this vertex we get
$$
p_a(n)= \sum_{g=0}^{N_1} \sum_{b\in\skrix_g} p_{a,b}(n-1-g)
\Q_r \{((0,\emptyset),g) \,|\, b\},
$$
implying that $S_a=\{n \, : \, n-1-g\in S_{a,b}$ for some $b\in\skrix_g\}$ and
for any $a \in \skrix_r$,
$$
\lim_{\heap{n\to\infty}{n\in S_a}} -\frac1{n}\log p_{a}(n)=I.
$$
Suppose for contradiction that $I>0$. Then, for
$a \in \skrix_r$ and
all $n\in S_a$
with $n\ge n_0$, we have $p_a(n)\le \exp(-nI/2)$. As $p_a(n)=0$ for
all $n\not\in S_a$, this implies that
$$\prob\big\{ |T|\ge n \, \big| \, X(\rho)=a\big\} \le
\frac{\exp(-nI/2)}{1-\exp(-I/2)}
\mbox{ for all } n\ge n_0.$$
But this probability is at least as large as the corresponding
probability for the restriction of $T$ to vertices whose type is in $\skrix_r$.
The latter is an irreducible, critical multitype Galton-Watson~tree, so
by the corollary in \cite[p.191]{AN72}
under the hypothesis of finite second moment this probability is
bounded below by a constant multiple of $1/n$, which is a
contradiction. Hence, $I=0$ and the result of the lemma follows
since by the weak irreducibility of $X$ we have that
$p_a(n)=0$ for all $n \geq n_0$ and $a \in \skrix_t$.
\end{proof}

\subsection{Proof of the upper bound in Theorem~\ref{general}}\label{s-upper}

Given a bounded 
function $\tilde{g}:\skrix\times \skrix^*\to\R$ we define the function
$$U_{\tilde{g}}(a)=\log \sum_{c\in \skrix^*} \Q\{ c\, | \, a\} e^{\tilde{g}(a,c)},$$
for $a\in\skrix$. We use $\tilde{g}$ to define a new multitype Galton-Watson tree
as follows:
\begin{itemize}
\item The type of the root $\rho$ is $a\in\skrix$ with probability
$$\mu_{\tilde{g}}(a) =\frac{e^{U_{\tilde{g}}(a)}\mu(a)}{\int e^{U_{\tilde{g}}(b)}
\mu(db)}.$$
\item for each vertex with type $a\in \skrix$ the offspring number
and types are given independently of everything else, by the
offspring law $\tilde{\Q}\{\,\cdot\, | \, a\}$ given by
\begin{align*}
\tilde{\Q}\big\{ c \,\big|\, a \big\} &  =
\exp\big(\tilde{g}\big(a,c\big)-U_{\tilde{g}}(a)\big) \Q\big\{c\,\big| \, a\big\}.
\end{align*}
\end{itemize}
We denote the transformed law by $\tilde{\prob}$ and make the simple observation
that $\tilde{\prob}$ is absolutely continuous with respect to
$\prob$, as for each finite $X\in\bar{\skrix}$, 
\begin{align}
\frac{d\tilde{\prob}}{d\prob}(X)& =
 \frac {e^{U_{\tilde{g}}(X(\rho))}}{\int e^{U_{\tilde{g}}(b)} \mu(db)}
\, \prod_{v\in V} \exp\Big[\tilde{g}(X(v),C(v))- U_{\tilde{g}}(X(v))\Big] \label{bddform}\\
& = \frac 1{\int e^{U_{\tilde{g}}(a)} \mu(da)}\prod_{v\in V}
\exp\Big[\tilde{g}(X(v),C(v))-\sum_{j=1}^{N(v)} U_{\tilde{g}}(X_j(v))
 \Big], \label{Jhatform}
\end{align}
recalling that $C(v)=\big(N(v),X_1(v),\ldots,X_N(v)\big)$.

We begin by establishing exponential tightness of the family of laws of $M_X$
on the space $\skrim(\skrix\times\skrix^*)$.

\begin{lemma}\label{tightness}
For every $A>0$ there exists a compact $K\subset\skrim(\skrix\times\skrix^*)$
with $$\limsup_{n\to\infty} \frac 1n \log \prob\big\{ M_X\not\in K \, \big| \,
|T|=n \big\} \le -A.$$
\end{lemma}

\begin{proof} Recall that $\Q\{ e^{\eta N} | \, a\}<\infty$ for all $\eta>0$.
Hence, given $l\in \mathbb{N}$, we may choose $k(l)\in \mathbb{N}$ so large that
$$\Q\big\{  \exp(l^2 N 1_{\{N>k(l)\}}) \, \big| \, a \big\} < 2 \mbox{
  for all } a\in\skrix.$$
Using the exponential Chebyshev inequality,
\begin{align*}
\prob\Big\{ & \int_{\{N>k(l)\}} N \, dM_X \ge \mbox{$\frac 1l$}, \, |T|=n \Big\}
\le  e^{-ln}\,
\me\Big\{ \exp\big( l^2 n  \int_{\{N>k(l)\}} N \, dM_X \big) ,\, |T|=n \Big\} \\
& = e^{-ln}
\me\Big\{ \prod_{v \in T} \exp\big( l^2 \1_{\{N(v)>k(l)\}} N(v)\big) ,\,
|T|=n \Big\} \\
& \le e^{-ln} \Big( \sup_{a\in\skrix}
\Q\big\{ \exp(l^2 N 1_{\{N>k(l)\}}) \, \big| \,  a\big\} \Big)^n \le e^{-n(l-\log 2)}.
\end{align*}
Now choose $M>A+\log 2$. Define the set
$$K= \Big\{ \nu\in\skrim(\skrix\times\skrix^*) \, : \,
\int _{\{N>k(l)\}} N \, d\nu <\mbox{$\frac 1l$, for all }l \ge M \Big\}.$$
As $\{N\le k(l)\}\subset\skrix\times\skrix^*$ is compact,
the set $K$ is pre-compact in the weak topology, by Prohorov's criterion.
Moreover, since $m(a,c)\le N$, it is easy to see by truncation that for every
weakly convergent sequence $\nu_n\to \nu$ with $\nu_n \in K$, we also have
$\lim_{n\to \infty} \int m(a,c) \nu_n(b,dc)=\int m(a,c) \nu(b,dc)$.
Hence, $K$ is even pre-compact in the stronger topology we are using
on the space $\skrim(\skrix\times\skrix^*)$.
As
$$\prob\big\{ M_X \not\in K \, \big| \, |T|=n \big\} \le
\frac{1}{\prob\{|T|=n\}} \frac{1}{1-e^{-1}} \exp(-n(M-\log 2)),$$
we can use Lemma~\ref{subexp} to infer that
$$\limsup_{n\to\infty} \frac 1n \log \prob\big\{ M_X\not\in K \, \big| \,
|T|=n \big\} \le -A,$$ as required for the proof.\end{proof}

Next we derive an upper bound in a variational formulation. Denote by
$\skric$ the space of bounded functions on  $\skrix\times \skrix^*$
and define for each $\nu \in \skrim(\skrix\times\skrix^*)$,
\begin{equation}\label{hatJ-def}
\widehat{J}(\nu) = \sup_{{g}\in \skric}\Big\{ \int
\Big[ {g}(b,c) -\sum_{j=1}^{n} U_{{g}}(a_j) \Big]
\, \nu(db\,, dc) \Big\},
\end{equation}
where $c=(n,a_1,\ldots,a_n)$.

\begin{lemma}\label{upper}
For each closed set $F\subset\skrim(\skrix\times\skrix^*)$,
\begin{align*}
\limsup_{n\to\infty} \frac 1n \log \prob\big\{ M_X \in F\, \big| \, |T|=n \big\} \le
-\inf_{\nu\in F} \widehat{J}(\nu).
\end{align*}
\end{lemma}

\begin{proof}
Fix $\tilde{g}\in\skric$ bounded by some $M>0$, then also
$\int e^{U_{\tilde{g}}(a)} \mu(da)\le e^M$. Define
$h:\skrix\times \skrix^*\to\R$ by
$h(b,c)=\tilde{g}(b,c)-\sum_{i=1}^n U_{\tilde{g}}(a_i)$, where
as usual $c=(n,a_1,\ldots, a_n)$, and observe that, by \eqref{Jhatform},
\begin{align*}
e^M & \, \ge \tilde{\prob}\{|T|=n\} \int e^{U_{\tilde{g}}(a)} \mu(da)
 = \me\Big\{ \prod_{v\in V} \exp\Big[\tilde{g}(X(v),C(v))-
\sum_{j=1}^{N(v)} U_{\tilde{g}}(X_j(v)) \Big] \1_{\{|T|=n\}} \Big\}\\
& = \me\big\{ e^{n \langle h, M_X \rangle} \, \1_{\{|T|=n\}} \big\}.
\end{align*}
Together with Lemma~\ref{subexp} this shows that
\begin{equation}\label{zero}
\limsup_{n\to \infty}\frac 1n \log \me\Big\{ e^{n \langle h, M_X \rangle} \,
\Big| \, |T|=n \Big\}\le 0.
\end{equation}
In view of \eqref{bddform} the same bound \eqref{zero} applies for
$h:\skrix\times \skrix^*\to\R$ of the form $h(b,c)=\tilde{g}(b,c)-U_{\tilde{g}}(b)$.

Now fix $\eps>0$, and let $\widehat{J}_\eps(\nu) = \min\{\widehat{J}(\nu),\eps^{-1}\}
-\eps$. Suppose first that 
$\nu\in F$ is shift-invariant. Then, for any
$\tilde{g}\in\skric$,
\begin{equation}\label{SIfact}
\int \sum_{j=1}^n U_{\tilde{g}}(a_j) \, \nu(db\, , dc) =
\sum_{(b,c)\in\skrix\times\skrix^*} \sum_{a\in\skrix} m(a,c)
\nu(b,c)U_{\tilde{g}}(a) = \sum_{a\in\skrix} U_{\tilde{g}}(a)\nu_1(a) =
\int U_{\tilde{g}}(b)\, \nu_1(db).
\end{equation}
Choose $\tilde{g}_\nu \in \skric$ such that
${h}_\nu (b,c)=\tilde{g}_\nu (b,c)-U_{\tilde{g}}(b)$
satisfies
$$ \langle {h}_\nu, \nu \rangle := \int {h}_\nu (b,c) \nu(db\, , dc)
=\int \big[ \tilde{g}_\nu(b,c) - \sum_{j=1}^n U_{\tilde{g}_\nu}(a_j)
\Big] \,\nu(db\, , dc) \geq \widehat{J}_\eps(\nu)  \,.$$
Since ${h}_\nu$ is bounded, the mapping $\langle {h}_\nu, \cdot \rangle$ is
continuous in $\skrim(\skrix\times\skrix^*)$. Hence there exists an open
neighbourhood $B_\nu$ of $\nu$ such that
$$\inf_{\mu \in B_{\nu}}\, \langle {h}_\nu, \mu \rangle \geq
\langle {h}_\nu, \nu \rangle - \eps \geq \widehat{J}_\eps (\nu)-\eps \, .$$
Using the exponential Chebyshev inequality and the remark following
\eqref{zero} we obtain that,
\begin{align}
\limsup_{n\to\infty} \frac 1n & \log \prob\big\{ M_X \in B_{\nu} \, \big|
\, |T|=n \big\} \nonumber\\
\le & \limsup_{n \to \infty}\frac{1}{n}
\log \me\big\{ e^{n \langle {h}_{\nu}, M_X \rangle} \,
\big| \, |T|=n \big\} - \widehat{J}_\eps(\nu) + \eps
\leq  - \inf_{\nu \in F} \widehat{J}_\eps(\nu) +  \eps.
\label{one}\end{align}
Now suppose that $\nu$ fails to be shift-invariant. Assume first that there exists $a\in \skrix$
such that
\begin{equation}\nu_1(a)< \sum_{(b,c)}m(a,c) \, \nu(b,c).
\label{notSI}\end{equation}
Recall that the mappings $\nu \mapsto \sum_{b,c} m(a,c) \, \nu(b, c)$
are continuous in our topology. Hence there exist $\delta>0$ and
a small open neighbourhood $B_\nu\subset\skrim(\skrix\times\skrix^*)$ such that
\begin{equation}\label{new-amir}
\tilde{\nu}_1(a)<\sum_{(b,c)}m(a,c) \, \tilde{\nu}(b,c)- \delta,
\mbox{ for all $\tilde{\nu}\in B_\nu$.}
\end{equation}
Let $\tilde{g}\in \skric$ be defined by
$\tilde{g}(b,c)=-(\delta \eps)^{-1}\1_a(b)$ and
$h(b,c)=\tilde{g}(b,c)-\sum_{j=1}^n U_{\tilde{g}}(a_j)$.
Note that
$U_{\tilde{g}}(b)=\tilde{g}(b,c)$ for all $b$ and
vanishes unless $b = a$. Hence, by \eqref{new-amir},
for every $\tilde{\nu}\in B_\nu$
we have that $\int h \, d\tilde{\nu}> \eps^{-1}$.
Then, using the exponential Chebyshev inequality and \eqref{zero},
\begin{align}
\limsup_{n\to\infty} \frac 1n & \log \prob\big\{ M_X \in B_{\nu} \, \big|
\, |T|=n \big\} \nonumber\\
\le & \limsup_{n \to \infty}\frac{1}{n}\log\me\big\{ e^{n \langle h, M_X \rangle}
\, \big| \, |T|=n \big\} - \eps^{-1} \le - \eps^{-1}
\leq  - \inf_{\nu \in F} \widehat{J}_\eps(\nu).
\label{two}\end{align}
In case the opposite inequality holds in \eqref{notSI} the same
argument leads to \eqref{two} if $\tilde{g}$ is defined as
$\tilde{g}(b,c)=(\delta \eps)^{-1}\1_a(b)$.

Now we use Lemma~\ref{tightness} to choose a compact set $K$ with
$$\limsup_{n\to\infty} \frac 1n \log \prob\big\{ M_X\not\in K \, \big| \,
|T|=n \big\} \le - \eps^{-1}.$$
The set $K \cap F$ is compact and hence it may be covered by finitely many
of the sets $B_{\nu_1},\ldots,B_{\nu_m}$, with $\nu_i \in F$ for $i=1,\ldots,m$.
Hence,
$$\prob\big\{ M_X \in F \, \big| \, |T|=n \big\}
\le \sum_{i=1}^m \prob\big\{ M_X \in B_{\nu_i} \, \big| \, |T|=n \big\}
+ \prob\big\{ M_X \not\in K \, \big| \, |T|=n \big\}.$$
Using \eqref{one} and \eqref{two} we obtain, for small enough $\eps>0$, that
$$\limsup_{n\to\infty} \frac 1n \log \prob\big\{ M_X \in F \, \big| \, |T|=n \big\}
\le \max_{i=1}^m \, \limsup_{n\to\infty} \frac 1n
\log \prob\big\{ M_X \in B_{\nu_i} \, \big| \, |T|=n \big\}
\le - \inf_{\nu \in F} \widehat{J}_\eps(\nu) +  \eps.$$
Taking $\eps \downarrow 0$ gives the required statement.
\end{proof}

We next show that the convex rate function~$J$ may replace
the function $\widehat{J}$ of \eqref{hatJ-def} in the upper bound
of Lemma~\ref{upper}.
\begin{lemma}\label{rate-lsc}
The function $J(\cdot)$ is convex and lower semicontinuous on
$\skrim(\skrix\times \skrix^*)$. Moreover,
$J(\nu) \leq \widehat{J}(\nu)$
for any $\nu \in \skrim(\skrix\times \skrix^*)$.
\end{lemma}

\begin{proof}
We start by proving the inequality $J(\nu) \leq \widehat{J}(\nu)$.
To this end, suppose first that $\nu\not\ll\nu_1\otimes\Q$.
Then, there exists $(a',c')\in\skrix\times\skrix^*$ with $\nu(a',c')>0$
and $\Q\{c'\,|\,a'\}=0$. Consequently, $U_{\tilde{g}}=0$
for $\tilde{g}(b,c)=K\1_{(a',c')}(b,c)$ and any $K$. Considering
such $\tilde{g}$ in \eqref{hatJ-def} with
$K \uparrow \infty$ we see that $\widehat{J}(\nu)=\infty$ in this case.

Suppose now that $\nu$ fails to be shift-invariant, in which case
there exists $a\in\skrix$ such that
$\nu_1(a) \neq \sum_{(b,c)\in \skrix\times \skrix^*} m(a,c) \nu(b,c)$.
Choose $\tilde{g}(b,c)=K \1_a(b)$, for which
$U_{\tilde{g}}(b)=K \1_a(b)$ and
$$\int \Big[ \tilde{g}(b,c)
-\sum_{j=1}^n U_{\tilde{g}}(a_j) \Big] \,\nu(db , dc)
= K \Big( \nu_1(a) - \int  m(a,c) \, \nu(db , dc) \Big)
\longrightarrow \infty,$$
for $|K|\uparrow\infty$, with the sign of $K$ chosen so that
the right hand side is positive.

Finally suppose that $\nu$ is shift-invariant and $\nu\ll\nu_1\otimes
\Q$. By the variational characterisation
of the relative entropy,
see e.g.~\cite[Lemma 6.2.13]{DZ98},
the definition of
$U_g$, Jensen's inequality, and \eqref{SIfact},
\begin{align}\label{ent-ident}
H(\nu\, \| \, \nu_1\otimes \Q)& = \sup_{g\in\skric} \Big\{ \int g\, d\nu
- \log \iint e^{g(a,c)}\, \Q\{dc\,|\, a\} \nu_1(da) \Big\}  \nonumber
\\
& = \sup_{g\in\skric} \Big\{ \int g\, d\nu
- \log \int e^{U_g(a)} \, \nu_1(da) \Big\}  \\
& \le \sup_{g\in\skric} \Big\{ \int g\, d\nu
- \int U_g(a)\, \nu_1(da) \Big\} = \widehat{J}(\nu).
\nonumber
\end{align}
If $\nu, \nu' \in \skrim(\skrix\times \skrix^*)$ are both
shift-invariant then $\nu_\lambda = \lambda \nu + (1-\lambda)\nu'$
is also shift-invariant for any $0<\lambda<1$. Moreover, $\nu \mapsto
\int m(a,c) \nu(b,dc)$ is continuous for each $a,b \in \skrix$, implying
that the set $\skris = \{ \nu : \nu $ is shift-invariant $\}$ is convex and
closed in the topology we use on $\skrim(\skrix\times \skrix^*)$.
Note that if $g \in \skric$, then so is $U_g$ and the mapping $\nu \mapsto
\int g d\nu - \log \int e^{U_g(a)} \nu_1(da)$ is continuous and
convex. Consequently, the identity \eqref{ent-ident} implies that
$\nu \mapsto H(\nu\, \| \, \nu_1\otimes \Q)$ is lower semicontinuous and
convex. For any $\alpha<\infty$, the level set $\{ \nu :
J(\nu) \leq \alpha \}$ is the intersection of the convex,
closed sets $\skris$ and
$\{ \nu : H(\nu\, \| \, \nu_1\otimes \Q) \leq \alpha \}$. Consequently,
$J(\cdot)$ is a convex rate function.
\end{proof}

\subsection{Proof of the lower bound in Theorem~\ref{general}}\label{s-lower}

Recall the definition of the multiplicity $m(a,c)$ of the symbol~$a$
in $c$ and of the matrix $A_{\tilde{g}}$ with index set $\skrix\times\skrix$
associated with the transformed multitype Galton-Watson tree,
$${A}_{\tilde{g}}(a,b)=\sum_{c\in \skrix^*}
\tilde{\Q}\{ c \, | \, b\} m(a,c),
\mbox{ for } a,b\in\skrix.$$
By our assumptions the matrix $A_{\tilde{g}}$
which has the same set of non-zero entries as $A$, is weakly irreducible.
Recall that, by the Perron-Frobenius theorem,
see e.g. \cite[Theorem 3.1.1]{DZ98}, the largest
eigenvalue $\varrho_{\tilde{g}}$
of the irreducible restriction of $A_{\tilde{g}}$ to $\skrix_r$
is real and positive, with strictly positive right and left eigenvectors.
Since $A_{\tilde{g}}$ is weakly irreducible,
the largest eigenvalue of
$A_{\tilde{g}}$ is also $\varrho_{\tilde{g}}$. Further, recall that
$A_{\tilde{g}} (a,b)=0$ whenever $b \in \skrix_t$ and $a \in \skrix_r$
or $b \leq a \in \skrix_t$, while
$\sum_{b \in \skrix_r} A_{\tilde{g}} (a,b)>0$ for any
$a \in \skrix_t$. Consequently, there exists
a unique right eigenvector $u_{\tilde{g}}\in\R^\skrix$ for the eigenvalue
$\varrho_{\tilde{g}}$ of $A_{\tilde{g}}$ having
strictly positive entries, which add up to one.
The next lemma guides the choice of $\tilde{g}$
associated with a
large deviations lower bound at
$\nu \in \skrim(\skrix\times \skrix^*)$ for which $J(\nu)<\infty$.

\begin{lemma}\label{shiftinv}
Suppose $\nu \in \skrim(\skrix\times \skrix^*)$
with $\nu_1$ strictly positive.
The following statements are equivalent.
\begin{itemize}
\item[(i)] $\nu$ is shift-invariant and $\nu \ll \nu_1 \otimes \Q$.
\item[(ii)] There exists a function $\tilde{g}:\skrix\times \skrix^*\to\R$
with $U_{\tilde{g}}=0$, such that $\varrho_{\tilde{g}}=1$ and the corresponding
Perron-Frobenius eigenvector~$u_{\tilde{g}}$ satisfies
$\nu(a,c)=\tilde{\Q}\{c \, | \, a\} u_{\tilde{g}}(a), \mbox{ for every }
(a,c)\in \skrix\times\skrix^*.$
\end{itemize}
Moreover, if (ii) holds, then $H(\nu\,\|\, \nu_1\otimes \Q)= \int
\tilde{g}(b,c) \, \nu(db\, , dc)$.
\end{lemma}

\begin{proof}
Suppose first that $\nu$ is shift-invariant and $\nu \ll \nu_1 \otimes \Q$.
Define $\tilde{g}$ by
\begin{equation}
\label{eq:gdef}
\tilde{g}(a,c)=\log\Big( \frac{\nu(a,c)}{\nu_1(a) \Q\{c\, | \,a\}}\Big)
\mbox{ when } \Q\{c\,|\,a\} > 0 \,,
\end{equation}
and otherwise $\tilde{g}(a,c)=0$. Then, for all $a\in\skrix$,
$$\sum_{c\in\skrix^*} \Q\{ c \, | \, a\}  e^{\tilde{g}(a,c)}= 1,$$
and hence $U_{\tilde{g}} (a)=0$. We infer that
\begin{equation}\label{tildeku}
\tilde{\Q}\{c \, | \, a\}=e^{\tilde{g}(a,c)}\Q\{c\, | \,a\}.
\end{equation}
Using this and the definition (\ref{eq:gdef})
of $\tilde{g}$ we see that
\begin{equation}\label{nuqu}
\nu(a,c)=e^{\tilde{g}(a,c)}\Q\{c\, | \,a\}\nu_1(a)=
\tilde{\Q}\{c \, | \, a\}\nu_1(a).
\end{equation}
To identify $\varrho_{\tilde{g}}$, by Perron-Frobenius theorem,
we only have to find
the eigenvalue corresponding to a strictly positive (right) eigenvector, which
turns out to be $\nu_1$. Indeed, for all $a\in\skrix$,
$$\sum_{b\in \skrix} A_{\tilde{g}}(a,b) \nu_1(b)
=  \sum_{(b,c)\in \skrix\times\skrix^*} \tilde{\Q}\{c \, | \, b\}
m(a,c) \nu_1(b) =  \sum_{(b,c)\in \skrix\times\skrix^*} \nu(b,c) m(a,c)= \nu_1(a),$$
using the shift-invariance of $\nu$ in the final step.
This shows that $\varrho_{\tilde{g}}=1$ and, by uniqueness of the eigenvector,
$\nu_1=u_{\tilde{g}}$. Hence (ii) follows from (\ref{nuqu}).

Conversely, fix $\tilde{g}$ for which $\varrho_{\tilde{g}}=1$ and (ii) holds.
Summing over $c\in\skrix^*$ in (ii) we have that $\nu_1 = u_{\tilde{g}}$
and hence $\nu \ll \nu_1 \otimes \Q$. Moreover, for all $a\in\skrix$,
$$\nu_1(a)= \sum_{b \in \skrix} A_{\tilde{g}}(a,b) \nu_1(b)=
\sum_{(b,c)\in \skrix\times\skrix^*} m(a,c)
\tilde{\Q}\{c \, | \, b\} \nu_1(b) =
\sum_{(b,c)\in \skrix\times\skrix^*} m(a,c) \nu(b,c),$$
hence $\nu$ is shift-invariant. Moreover, 
using $\nu(a,c)=\tilde{\Q}\{c \, | \, a\}\nu_1(a)$ and the definition of
$\tilde{\Q}$, we get
$$H(\nu\, \| \, \nu_1\otimes \Q)  =
\sum_{(a,c)\in\skrix\times\skrix^*} \nu(a,c)
\log \frac{\tilde{\Q}\{c\, | \,a\}}{\Q\{c\, | \,a\}} =
\int \tilde{g}(a,c) \, \nu(da, dc),$$
which completes the proof.
\end{proof}

The next lemma is key to the proof of
the lower bound in Theorem~\ref{general}. It allows us
to focus on 
those shift-invariant $\nu \in \skrim(\skrix\times\skrix^*)$
with strictly positive first marginal, for which $\tilde{g}$ of Lemma
\ref{shiftinv} is bounded above. 
If $\nu \in\skrim(\skrix\times\skrix^*)$ and $a \in \skrix$ we write
$\nu(\,\cdot\,|\, a)= \nu(\,\cdot\,, a)/\nu_1(a)$.

\begin{lemma}\label{apx} 
Suppose $O$ is an open subset of $\skrim(\skrix\times\skrix^*)$
and $\nu\in O$ with $J(\nu)<\infty$. Then, for any $\delta>0$, there exists
$\tilde{\nu} \in O$ with $J(\tilde{\nu}) \leq J(\nu) + \delta$, such that
$\tilde{\nu}_1$ is strictly positive and
$\tilde{\nu}(c\,|\,a) \leq \Q\{c\,|\,a\}/y$ for some $y>0$ and all $(a,c) \in
\skrix\times\skrix^*$.
\end{lemma}

\begin{proof} Recall our assumption that $X$ is
weakly
irreducible and critical. This implies the existence
of a strictly positive probability vector $u_0$ on $\skrix$ such that
$\nu^*(a,c)=\Q\{c\,|\,a\} u_0(a) \in \skrim(\skrix\times\skrix^*)$
is shift-invariant with $\nu^*_1(a)=u_0(a)$ and $J(\nu^*)=0$.
Fixing $\nu \in O$ with $J(\nu)<\infty$, we have for each $0< \eps < 1$ that
$\nu_\eps = (1-\eps) \nu + \eps \nu^*$ is
shift-invariant in $\skrim(\skrix\times\skrix^*)$
with $(\nu_\eps)_1$ strictly positive and
$\nu_\eps (c\,|\,a)=0$ exactly for 
those values $(a,c)\in\skrix\times\skrix^*$ where $\Q\{c\,|\,a\}=0$.
By convexity of $J(\cdot)$ we know that $J(\nu_\eps) \leq (1-\eps) J(\nu)$.
Further, $\int f d \nu_\eps \to \int f d\nu$ as $\eps \downarrow 0$,
for any $f:\skrix\times\skrix^*\to\R$ which is either bounded
or satisfies $f(b,c)=m(a,c) \1_{b_0}(b)$ for some
$a,b_0 \in \skrix$. As $O$ is open in
$\skrim(\skrix\times\skrix^*)$, it follows that
$\nu_\eps \in O$ for all $\eps>0$ small enough.

In view of the above, we may and shall assume hereafter that $\nu_1$ is
strictly positive and $\nu (c\,|\,a)=0$ exactly for those values
$(a,c)\in\skrix\times\skrix^*$  where $\Q\{c\,|\,a\}=0$. In particular,
the matrix $ A_{0,0}$ given by
$$A_{0,0}(a,b)=\sum_{c\in\skrix^*} m(a,c) \nu(c\,|\,b),
\mbox{ for } a,b\in\skrix,$$
has nonnegative entries and is
weakly
irreducible.
Its Perron-Frobenius eigenvalue, denoted $\varrho(A_{0,0})$,
equals $1$, and the corresponding right eigenvector $u_{0,0}$
equals $\nu_1$ and hence is a strictly positive probability vector on
$\skrix$. The corresponding left eigenvector $v_{0,0}$ is a probability
vector which is strictly positive on $\skrix_r$.
Clearly, for each $b \in \skrix_r$ there exists
$c_1=c_1(b)$ such that $\Q\{c_1\,|\,b\}>0$, hence also $\nu(c_1\,|\,b)>0$.
Recall that for $b \in \skrix_t$ we have $\Q\{c\,|\,b\}>0$ (and hence
$\nu(c\,|\,b)>0$) for only finitely many $c \in \skrix^*$. Consequently,
$\nu(c\,|\,b) \leq \Q\{c\,|\,b\}/y$ for some $y>0$ and all $c \in \skrix^*$,
$b \in \skrix_t$. The proof of the lemma is complete if the same
applies for all $b \in \skrix_r$. Assuming hereafter that this is
not the case, with $\sum_{a \in \skrix_t} m(a,c)$ uniformly bounded
under $\Q$, there must exist
$b_0 \in \skrix_r$ and $c_2=c_2(b_0) \in \skrix^*$
such that $\Q\{c_2\,|\,b_0\}>0$ (and hence
also $\nu(c_2\,|\,b_0)>0$), with $\sum_{a \in \skrix_r} m(a,c_2)$
large enough to guarantee that
$\sum_{a \in \skrix_r} v_{0,0} (a) (m(a,c_2)-m(a,c_1(b_0))) >0$.
Let $c_1(b)$ be arbitrary for $b \in \skrix_t$, and $c_2=c_1(b)$
for all $b \neq b_0$.

Using these $c_1$ and $c_2$
we next construct probability measures
$\nu_{x,y}(\,\cdot\,|\,b)$ on $\skrix^*$ for $0< y < y_0$ and $|x| < 1/2$,
such that for each $b \in \skrix$ and $c \in \skrix^*$ we have
\begin{itemize}
\item $\nu_{x,y}(c\,|\,b) \leq \Q\{c\,|\,b\}/y$,
\item $\nu_{x,y}(c\,|\,b) \to \nu_{0,0} (c\,|\,b) = \nu(c\,|\,b)$
as $x \to 0$ and $y \downarrow 0$,
\item $\nu_{x,y}(c\,|\,b)=0$ if and only if $\nu(c\,|\,b)=0$.
\end{itemize}
Further,
\begin{equation}\label{usc-xy}
\limsup_{\heap{x \to 0}{y \downarrow 0}} H\big(\nu_{x,y} (\cdot\,|\,b)\, \big\| \,
\Q\{\cdot\,|\,b\} \big)
\leq H\big(\nu_{0,0} (\, \cdot\,|\,b) \,\big\|\, \Q\{\cdot\,|\,b\} \big) \,,
\end{equation}
and $A_{x,y}(a,b)=\sum_c m(a,c) \nu_{x,y} (c\,|\,b)  \to  A_{0,0} (a,b)$
for any $a,b \in \skrix$.
Note that $A_{x,y}(a,b)=0$ if and only if
$A_{0,0}(a,b)=0$, so with $A_{0,0}$
weakly irreducible, the same applies to $A_{x,y}$.
The function $f(x,y)=\varrho(A_{x,y})$ is thus continuous in this range of
$(x,y)$, as is also the strictly positive Perron-Frobenius right eigenvector
$u_{x,y}$ of $A_{x,y}$, normalized to be a probability vector on $\skrix$. Our
construction is such that $A_{x,0} = A_{0,0} + x B$
where $B(a,b)=\nu(c_2\,|\,b) \nu(c_1\,|\,b) (m(a,c_2)-m(a,c_1))$.
Therefore, $f(x,0)$ is continuously differentiable at $x=0$ with
$$
\frac{\partial f}{\partial x}(0,0) =
\frac{\sum_{a,b} v_{0,0}(a) B(a,b) u_{0,0} (b)}
{\sum_a v_{0,0}(a) u_{0,0}(a)} > 0 \,.
$$
By the implicit function theorem, there exist $x(y) \to 0$
as $y \downarrow 0$ such that $f(x(y),y)=f(0,0)=1$ for all $y>0$
small enough. It follows that
$\nu_{x,y}(b,c)=\nu_{x,y}(c\,|\,b) u_{x,y}(b)$ defines a
shift-invariant probability measure $\nu_{x,y}\in\skrim(\skrix\times\skrix^*)$
for $x=x(y)$ and all $y>0$ small enough. Moreover,
$$\int m(a,c) \nu_{x(y),y} (b,dc)=A_{x(y),y}(a,b) u_{x(y),y}(b)
\to A_{0,0}(a,b) u_{0,0}(b) =  \int m(a,c) \nu(b,dc) \,,$$
for each $a,b \in \skrix$ and $y \downarrow 0$, implying the convergence of
$\nu_{x(y),y}$ to $\nu$ in the topology of $\skrim(\skrix \times \skrix^*)$,
and by \eqref{usc-xy} and shift-invariance, also 
\begin{align*}
\limsup_{y \downarrow 0} J(\nu_{x(y),y}) & =
\limsup_{y \downarrow 0} \sum_{b\in\skrix} u_{x(y),y}(b)
H\big(\nu_{x(y),y} (\,\cdot\,|\,b) \,\big\|\, \Q\{\,\cdot\,|\,b\} \big)\\
& \leq \sum_{b\in\skrix} u_{0,0} (b)
H\big(\nu_{0,0} (\,\cdot\,|\,b) \,\big\|\, \Q\{\,\cdot\,|\,b\} \big) = J(\nu) \,,
\end{align*}
which completes the proof of the lemma 
subject to the construction of $\nu_{x,y}(\,\cdot\,\,|\,b)$.

We now turn to this construction. For any $|x| < 1/2$ we define the probability
measure
$$\nu_{x,0}(c\,|\,b)=\nu(c\,|\,b)+x \nu(c_2\,|\,b)\nu(c_1\,|\,b)
( \1_{\{c=c_2\}} - \1_{\{c=c_1\}} ) \,. $$
In particular, $\nu_{x,0}(c\,|\,b)=0$ exactly where $\nu(c\,|\,b)=0$ and
$A_{x,0}=A_{0,0}+x B$ as stated. Let
$y_0 = \Q\{c_2\,|\,b_0\} \min_{b \in \skrix_r} \Q\{c_1\,|\,b\} >0$ further
reducing $y_0$ as needed 
to ensure that $\nu(c\,|\,b) \leq \Q\{c\,|\,b\}/y_0$ for any $c \in \skrix^*$
and $b \in \skrix_t$. For any $0<y<y_0$
define the probability measures 
$\nu_{x,y}(\,\cdot\,|\, b \,)$ by
\begin{align*}
\nu_{x,y}(c\,|\,b)&=\min( \nu_{x,0}(c\,|\,b), \Q\{c\,|\,b\}/y )
\mbox{ for } c \neq c_1\,, \nonumber \\
\nu_{x,y}(c_1\,|\,b)&=\nu_{x,0}(c_1\,|\,b) +
\sum_{c \neq c_1} (\nu (c\,|\,b) -\Q\{c\,|\,b\}/y)_+ \, ,
\end{align*}
with $_+$ indicating the positive part. Our choice of $y_0$ results in
$\nu_{x,y}(\,\cdot\,|\,b)=\nu(\,\cdot\,|\,b)$ whenever $b \in \skrix_t$
and further guarantees that
$$\nu_{x,y}(c_2\,|\,b_0) = \nu_{x,0}(c_2\,|\,b_0) \leq \Q\{c_2\,|\,b_0\}/y$$
and $\nu_{x,y}(c_1\,|\,b) \leq 1 \leq \Q\{c_1\,|\,b\}/y$
for all $b \in \skrix_r$, $|x|<1/2$ and $0<y<y_0$.
Hence we have as stated that
$\nu_{x,y}(c\,|\,b) \leq \Q\{c\,|\,b\}/y$ for all $c \in \skrix^*$,
and $\nu_{x,y}(c\,|\,b)=0$ if and only if $\nu(c\,|\,b)=0$.
Moreover, $A_{x,y}=A_{x,0}+E_y$,
for 
$$E_y(a,b)=\sum_{c\in\skrix^*}
(m(a,c_1)-m(a,c)) \big(\nu(c\,|\,b)-\Q\{c\,|\,b\}/y \big)_+,$$
in particular, $E_y(a,b)=0$ for $b \in \skrix_t$.
Writing $n(c)=n$ if $c\in\skrix^n$.
Recall that $\sum_c n(c) \nu(c\,|\,b) = \sum_a A_{0,0}(a,b) < \infty$
for all $b \in \skrix$, so by dominated convergence
$$|E_y(a,b)| \leq \sum_{c\in\skrix^*} (n(c_1)+n(c)) \nu(c\,|\,b)
\1_{\{\nu(c\,|\,b) > \Q\{c\,|\,b\}/y\}}
\quad \underset{y\downarrow 0}{\longrightarrow} \quad 0 \,,$$
and consequently, as stated, each entry of $A_{x,y}$ is continuous in
$(x,y) \in (-1/2,1/2) \times [0,y_0)$. By the same argument,
$\sum_{c \neq c_1} (\nu (c\,|\,b)-\Q\{c\,|\,b\}/y)_+ \to 0$ as $y \downarrow 0$,
implying the pointwise convergence $\nu_{x,y}(c\,|\,b) \to \nu (c\,|\,b)$
for each $(b,c) \in \skrix \times \skrix^*$. Turning to
\eqref{usc-xy}, note that it suffices to consider only $b \in \skrix_r$.
Recall that
for any $q>0$ the function $z \log(z/q)$ increases in $z \in [q,1]$, and if
$\nu_{x,y}(c\,|\,b) \neq \nu_{0,0}(c\,|\,b)$ and $c \neq c_1$, $c \neq c_2$,
then necessarily $0<\Q\{c\,|\,b\} \leq \nu_{x,y}(c\,|\,b) < \nu_{0,0}(c\,|\,b) < 1$.
Consequently,
$$\sum_{\heap{c \neq c_1}{c \neq c_2}}
\nu_{x,y} (c\,|\,b) \log \frac{\nu_{x,y} (c\,|\,b)}{\Q\{c\,|\,b\}}
\leq
\sum_{\heap{c \neq c_1}{c \neq c_2}}
\nu_{0,0} (c\,|\,b) \log \frac{\nu_{0,0} (c\,|\,b)}{\Q\{c\,|\,b\}} \,,
$$
yielding \eqref{usc-xy} since
$\nu_{x,y} (c_i\,|\,b) \to  \nu_{0,0} (c_i\,|\,b)$
and $\Q\{c_i\,|\,b\}>0$
for $i=1,2$ and $b \in \skrix_r$.
\end{proof}

Using Lemma \ref{apx} we now establish the lower bound in
Theorem \ref{general}.
\begin{lemma}
For each open set $O\subset\skrim(\skrix\times\skrix^*)$,
$$\liminf_{n\to\infty}
\frac 1n \log \prob\big\{ M_X \in O \,  \big| \, |T|=n \big\}
\ge -\inf_{\nu\in O} J(\nu). $$
\end{lemma}

\begin{proof}
Suppose that $\nu$ is an approximate minimizer on the right hand side.
We can assume without loss of generality that $J(\nu)<\infty$,
hence $\nu$ is shift-invariant with $\nu \ll \nu_1 \otimes \Q$.
By Lemma \ref{apx} we may and shall assume in addition that
$\nu_1$ is strictly positive and the function
$\tilde{g}$ 
associated to $\nu$ via (\ref{eq:gdef}) is bounded from above.
Recall from Lemma~\ref{shiftinv} that $\varrho_{\tilde{g}}=1$, and the
corresponding Perron-Frobenius eigenvector~$u_{\tilde{g}}$ satisfies
$$\nu(a,c)=\tilde{\Q}\{c \, | \, a\} \,u_{\tilde{g}}(a), \mbox{ for every }
(a,c)\in \skrix\times\skrix^*,$$
and further that
$H(\nu\, \| \, \nu_1\otimes \Q) = \int \tilde{g}(b,c)\, \nu(db , dc).$
It thus suffices to show that
$$\liminf_{n\to\infty} \frac 1n \log \prob\big\{ M_X \in O \,
\big| \, |T|=n \big\} \ge -\int \tilde{g}(b,c) \, \nu(db, dc).$$
Since $\tilde{g}$ is bounded above, fixing $\eps>0$ we can
choose an open set $\widetilde{O} \subset O$ such that $\nu \in \widetilde{O}$
and $\langle \tilde{g}, \mu\rangle \le  \langle \tilde{g}, \nu\rangle + \eps$
for all $\mu\in \widetilde{O}$.
We use the transformed probability measures $\tilde{\prob}$ and the formula
\eqref{bddform} for their density, to get
\begin{align*}
\prob\big\{ M_X & \in O , \, |T|=n \big\}
\ge \tilde{\me}\bigg\{ \frac{d\prob}{d\tilde{\prob}}(T) \1_{\{M_X\in
\widetilde{O}\}}
\1_{\{|T|=n\}} \bigg\} \\
& = \tilde{\me}\bigg\{ \prod_{v\in V}
\exp\Big(-\tilde{g}(X(v),C(v))\Big) \1_{\{M_X\in
\widetilde{O} \}} \1_{\{|T|=n\}} \bigg\} \\
& \ge \exp\big(-n \langle \tilde{g}, \nu\big\rangle
- n\eps \big) \times \tilde{\prob}\Big\{ M_X\in
\widetilde{O} ,\, |T|=n\Big\}.
\end{align*}
Dividing by $\prob\{|T|=n\}$ and recalling Lemma~\ref{subexp} gives
\begin{align*}
\liminf_{n\to\infty} & \frac 1n \log \prob\big\{ M_X \in O \, \big| \,
 |T|=n \big\}\\
& \ge -n \langle \tilde{g}, \nu\big\rangle- n\eps
+ \liminf_{n\to\infty} \frac 1n \log \tilde{\prob}\Big\{ M_X\in
\widetilde{O} \Big| \, |T|=n\Big\}.
\end{align*}
The result follows once we show that
\begin{equation}\label{toshow}
\limsup_{n\to\infty} \frac 1n \log \tilde{\prob}\Big\{
M_X\notin
\widetilde{O} \Big| \, |T|=n\Big\}<0.
\end{equation}
We use the upper bound (but now with the law $\prob$ replaced
by $\tilde{\prob}$) to establish (\ref{toshow}). Indeed, since $\tilde{g}$ is
bounded from above, we have $\tilde{\Q} \{e^{\eta N} \,| \, a\} < \infty$ for all
$a \in \skrix$ and $\eta>0$. So, denoting
$$\tilde{J}(\nu)=\left\{ \begin{array}{ll} H(\nu \,\|\,  \nu_1\otimes\tilde{\Q}) &
\mbox{ if $\nu$ is shift-invariant,} \\ \infty & \mbox{ otherwise,}
\end{array} \right.$$
the upper bound gives
$$\limsup_{n\to\infty} \frac 1n
\log \tilde{\prob}\Big\{ M_X\notin  \widetilde{O}\, \Big| \, |T|=n\Big\}
\le - \inf_{\tilde{\nu}\in K} \tilde{J}(\tilde{\nu}), $$
where $K\subset \widetilde{O}^c$
is a compact subset of $\skrim(\skrix\times\skrix^*)$. It suffices to show that
the infimum is positive. Suppose, for contradiction, that there exists a
sequence $\tilde{\nu}_n$ with $\tilde{J}(\tilde{\nu}_n)\downarrow 0$. By compactness
of $K$ and lower semicontinuity of $\nu\mapsto \tilde{J}(\nu)$, we can
extract a limit point $\tilde{\nu}\in K$ with $\tilde{J}(\tilde{\nu})=0$, and
hence $\tilde{\nu}$ is shift-invariant and $H(\tilde{\nu}\, \| \,
\tilde{\nu}_1 \otimes \tilde{\Q})=0$. This implies that
$\tilde{\nu}(a,c)=\tilde{\Q}\{c\, | \, a\}
\tilde{\nu}_1(a)$, for every $(a,c)\in \skrix\times\skrix^*$. Then, using
shift-invariance of $\tilde{\nu}$, for any $b\in\skrix$,
$$\sum_{(a,c)\in\skrix\times\skrix^*}
\tilde{\Q}\{c\,|\,a\} m(b,c)\tilde{\nu}_1(a)
= \sum_{(a,c)\in\skrix\times\skrix^*}
\tilde{\nu}(a,c) m(b,c)=\tilde{\nu}_1(b).$$
By the uniqueness of the Perron-Frobenius eigenvector we infer that
$\tilde{\nu}_1=u_{\tilde{g}}=\nu_1$ and this implies $\tilde{\nu}=\nu$, which
contradicts $\tilde{\nu}\in K$.
\end{proof}

We complete the proof of Theorem~\ref{general} by noting that the rate
function $J$ has compact level sets, 
i.e. is a \emph{good} rate function. This follows from abstract considerations
as stated, e.g., in \cite[Theorem 1.2.18]{DZ98}.

\subsection{Proof of Theorem~\ref{main}}\label{s-main}

Note that $X$ is an irreducible, critical multitype Galton-Watson
tree with offspring law
$$\Q\{c\,|\,b\}=p(n)\prod_{i=1}^{n} Q\{a_i\,|\,b\},
\mbox{  for $c=(n,a_1,\ldots,a_n)$,}$$
such that all exponential moments are finite. We derive Theorem~\ref{main}
from Theorem~\ref{general} by applying the contraction principle to
the continuous linear mapping
$F:\skrim(\skrix \times \skrix^*) \to \R^{\skrix \times \skrix}$,
defined by
$$F(\nu)(a,b)=\sum_{c\in\skrix^*} m(b,c) \nu(a,c)
\mbox{ for all } \nu \in \skrim(\skrix \times \skrix^*) \mbox{ and }
a,b\in\skrix.$$
Indeed, Theorem \ref{general}
implies the large deviation principle for $F(M_X)$ conditioned on
$\{ |T|=n \}$ with the good rate function $I(\mu)=\inf\{J(\nu) :
F(\nu)=\mu \}$, see for example \cite[Theorem 4.2.1]{DZ98}.
Convexity of $I$ follows easily from the linearity of $F$ and
convexity of $J$. It is easy to see that on  $\{ |T| = n \}$
we have $L_X=\frac{n}{n-1} F(M_X)$.
It follows that conditioned on  $\{ |T| = n \}$ the random variables $L_X$
are exponentially equivalent to $F(M_X)$, hence $L_X$ satisfy
the same large deviation principle as $F(M_X)$, see \cite[Theorem 4.2.13]{DZ98}.
Without loss of generality we restrict the space for the
large deviation principle of $L_X$ to the set of all probability vectors
on $\skrix \times \skrix$,
see \cite[Lemma 4.1.5(b)]{DZ98}.

Turning to the proof of \eqref{Idef}, recall that $\nu$ is shift-invariant
if and only if $\sum_a F(\nu)(a,b)=\nu_1(b)$ for all $b \in \skrix$. Hence,
if also $F(\nu)=\mu$, then necessarily $\nu_1=\mu_2$ and consequently,
$$I(\mu) = \inf\big\{ H(\nu \, \| \nu_1\otimes\Q)\, : \, F(\nu)=\mu,
\nu_1=\mu_2 \big\} \;.$$
Note that $\nu_1(a)=0$ yields $\sum_b F(\nu)(a,b)=0$. Hence if
$\mu_1(a)>0=\mu_2(a)$ for some $a \in \skrix$ then
$\{\nu : F(\nu)=\mu, \nu_1 = \mu_2\}$ is an empty set, and therefore
$I(\mu)=\infty$. Assuming hereafter that $\mu_1 \ll \mu_2$, it is
not hard to check that
\begin{equation}\label{I-abs}
I(\mu) = \sum_{a\in\skrix} \mu_2(a) \,
\widetilde{I} \Big(\frac{\mu(a,\cdot)}{\mu_2(a)},\Q\{\,\cdot\,|\,a\}\Big) \;,
\end{equation}
where for $\phi:\skrix \to \reals_+$ and $q \in \skrim(\skrix^*)$,
\begin{equation}\label{Ia-dec}
\widetilde{I}(\phi,q) = \inf\Big\{ H(\widetilde{\nu} \,\|\,  q)\, : \,
\widetilde{\nu} \in \skrim(\skrix^*), \;\;
\phi(b)=\sum_{c\in\skrix^*} m(b,c) \, \widetilde{\nu} (c) \mbox{ for all }
b \in \skrix \Big\} \;.
\end{equation}
Suppose now that $q(c)=p(n) \prod_{i=1}^n \widehat{q} (a_i)$
for all $c=(n,a_1,\ldots,a_n)$, where $\widehat{q}(\cdot)$ is a probability
vector on $\skrix$ and $p(\,\cdot\,)$ a probability measure with mean one
on the nonnegative integers, whose exponential moments are all finite.
With $z=\sum_b \phi(b)$ we show next that,
\begin{equation}\label{Ia-iden}
\widetilde{I}(\phi,q) = z H\big(\phi/z \,\|\,  \widehat{q}\big) + I_p(z) \;.
\end{equation}
Once this is done, we combine \eqref{Ia-iden} for $\widehat{q}(\cdot)=
Q\{\,\cdot\,|\,a\}$ and $z=\mu_1(a)/\mu_2(a)$ with the representation \eqref{I-abs}
of $I(\mu)$, which directly yields the formula \eqref{Idef}, thus completing the
proof of the theorem.

To prove \eqref{Ia-iden}, suppose first that $z=0$, i.e.
$\phi(b)=0$ for all $b \in \skrix$.
In this case, $\widetilde{\nu}((0,\emptyset))=1$ is the only possible
measure in \eqref{Ia-dec}, leading to
$\widetilde{I}(\phi,q) =-\log q((0,\emptyset))= -\log p(0)$, whereas
it follows from \eqref{Ip-def} that $I_p(0)=-\log p(0)$ establishing
\eqref{Ia-iden} for such $\phi(\cdot)$. Assume hereafter that $z>0$.
Now the possible measures $\widetilde{\nu}(\cdot)$ in \eqref{Ia-dec} are
of the form $\widetilde{\nu}(c)=s(n) v_n(a_1,\ldots,a_n)$
for $c=(n,a_1,\ldots,a_n)$, with $v_0=1$, where $s(\cdot)$ is
a probability measure on the nonnegative integers whose mean
is $z$, and $v_n(\,\cdot\,)$, $n \geq 1$, are probability measures on $\skrix^n$
with marginals $v_{n,i}(\,\cdot\,)$ such that
\begin{equation}\label{Ia-cons}
\phi(b)=\sum_{n=1}^\infty s(n) \sum_{i=1}^n v_{n,i} (b)
\quad  \mbox{ for all } b \in\skrix \;.
\end{equation}
By the assumed structure of $q(\,\cdot\,)$ we have
for such $\widetilde{\nu}(\,\cdot\,)$ that
$$ H(\widetilde{\nu} \,\|\,  q)  =
\sum_{n=1}^\infty s(n) H ( v_n \, \| \, \widehat{q}^n ) +
H( s \, \| \, p)\; ,$$
where $\widehat{q}^n$ denotes the product measure on $\skrix^n$ with
equal marginals $\widehat{q}$. Recall that
$$\sum_{n=1}^\infty s(n) H \big( v_n \, \big\| \, \widehat{q}^n \big) \geq
\sum_{n=1}^\infty s(n) \sum_{i=1}^n H\big(v_{n,i}\,\big\|\, \widehat{q}\big)
 \geq z H\,\Big( z^{-1}
\sum_{n=1}^\infty s(n) \sum_{i=1}^n v_{n,i}\,\Big\|\, \widehat{q}\Big) \,,$$
with equality whenever $v_n = \prod_{i=1}^n v_{n,i}$ and
$v_{n,i}$ are independent of $n$ and $i$ (see \cite[Lemma 7.3.25]{DZ98}
for the first inequality, with the second inequality following by
convexity of $H(\cdot\,\|\, \widehat{q})$ and the fact
that $\sum_n s(n) n =z$). So, in view of \eqref{Ia-cons},
\begin{equation}\label{Ia-temp}
H(\widetilde{\nu} \,\|\,  q)  \geq z H(\phi/z \,\|\,  \widehat{q}) +
H( s \, \| \, p) \,,
\end{equation}
with equality when
$v_n = (z^{-1} \phi)^n$ for all $n \geq 1$. Recall that with
all exponential moments of $p(\cdot)$ finite,
$I_p(z) = \inf \{ H( s \, \| \, p) : s(\cdot)$
a probability measure on $\{0,1,\ldots\}$ and
$\sum_n s(n) n = z \}$ (see \cite[(2.1.27)]{DZ98} for a similar identity).
Combining this with \eqref{Ia-temp} leads to \eqref{Ia-iden} and
completes our proof.

\subsection{Proof of Theorem~\ref{process}}\label{s-process}

In the first step we
extend the result of
Theorem~\ref{general}
to $k$-generation empirical offspring measures,
for each $k \geq 2$, in case
$\Q$ is irreducible and
the offspring size is bounded by some
non-random $N_0<\infty$.

For each $k\ge 0$, let $\skrix(k)$ be the
{\it finite}
set of typed trees with height at most $k$ and maximal degree $N_0+1$,
equipped with the discrete topology
(in particular, $\skrix(0)=\skrix$).
Let $\pi_k:\bar{\skrix}\to\skrix(k)$ be the canonical projection
obtained by removing all vertices in generations exceeding $k$
and
$\pi_{k,l}:\skrix(k) \to \skrix(l)$, $k \geq l$, the projections
obtained by removing all vertices in generations exceeding $l$.

If $X$ is a finite typed tree and $v$ is a vertex in this tree, we denote
by $X^v$ the subtree rooted in $v$
and let the $k$-generation empirical offspring
measures~$M_X^k$ associated to $X$ be defined as
$$
M_X^k(b)=\frac 1{|T|} \sum_{v\in V} \delta_{\pi_k(X^v)} (b),
\mbox{ for all }b\in{\skrix(k)}
$$
(for example $M_X^1(b)=M_X(a,c)$ where $b \in \skrix(1)$ has root of type
$a$ with $n$ children of types $a_1,\ldots,a_n$ and $c=(n,a_1,\ldots,a_n)$).
Given $a\in \skrix(k-1)$ and $b\in \skrix(k)$ we write $m_k(a,b)$ for
the number of children $v$ of the root in $b$ such that $b^v=a$.
A measure $\mu$ on $\skrix(k)$ is called \emph{shift-invariant} if
\begin{equation}\label{disp0}
\mu\circ\pi_{k,k-1}^{-1}(a) = \sum_{b\in \skrix(k)} m_k(a,b) \mu(b),
\mbox{ for all } a \in \skrix(k-1).
\end{equation}
We equip the space
$\skrim(\skrix(k))$
of probability measures on $\skrix(k)$ with the
smallest topology which makes the functionals $\mu \mapsto \int f d\mu$
continuous for each bounded $f:\skrix(k)\to\R$
(since the maximal degree is bounded in $\skrix(k)$, it follows that
$\mu \mapsto
\int m_k(a,x) d\mu(x)$ is also continuous for each $a \in \skrix(k-1)$).

Define $\mu\circ\pi_{k,k-1}^{-1} \otimes_1 \Q$ as the measure on $\skrix(k)$
obtained by providing
children for each vertex of the $k-1$ generation, independently according
to the transition mechanism $\Q$, and define the function
$$J_k(\mu)=\left\{ \begin{array}{ll}
H\big(\mu \, \big\| \, \mu\circ\pi_{k,k-1}^{-1}\otimes_1 \Q \big) &
\mbox{ if $\mu$ is shift-invariant,} \\ \infty & \mbox{ otherwise,}
\end{array} \right.$$
on $\skrim(\skrix(k))$. Note that $J_1(\cdot)$ coincides with
the good rate function $J(\cdot)$ of Theorem \ref{general}.

\begin{lemma}\label{kstep}
Suppose that $X$ is an irreducible, critical multitype Galton-Watson tree
with uniformly bounded offspring sizes,
conditioned to have exactly $n$ vertices. Then, for $n\to\infty$,
the $k$-generation empirical offspring measure~$M^k_X$ satisfies a large deviation
principle in $\skrim(\skrix(k))$ with speed $n$ and
convex,
good rate function $J_k(\cdot)$.
\end{lemma}

\begin{proof}
For $l \geq 0$ let $\skrix \{ l \} \subset \skrix (l)$ be the support of
$\pi_l(X)$ for a multitype Galton-Watson tree $X$ corresponding to the
transition mechanism $\Q$ starting at any strictly positive measure for $X(\rho)$.
Let $\skrix_m\{ l \}$ be the partition of $\skrix\{l\}$
according to the height $m=0,1,\ldots,l$ of the tree.
Let
$$\skrii: \skrix \{ k \} \to \skrix \{ k-1 \} \times \skrix \{ k-1 \}^*
\mbox{ given by } \left\{\begin{array}{lcl}
\skrii_1 (b)& = &\pi_{k,k-1}(b) \in \skrix \{ k-1 \},\\
\skrii_2 (b)& = & (n,b^{v_1},\ldots,b^{v_n}) \in \skrix \{ k-1 \}^*,
\end{array} \right.$$
where $v_1,\ldots,v_n$ are the vertices in the first generation of
$b\in\skrix \{ k \}$ ordered from left to right. 

To prove Lemma~\ref{kstep} we intend to apply Theorem~\ref{general} to a
multitype Galton-Watson tree $\widetilde{X}$ on the enlarged finite
type space $\skrix \{ k-1 \}$. We mark the objects related to this
new tree by\ \ $\widetilde{}$.

The process $\widetilde{X}$ is constructed by choosing
$\widetilde{X} (\rho)$ using the law of $\pi_{k-1} (X)$,
and the offspring number and types of a vertex $v$ as
$\widetilde{C}(v)=\skrii_2 (b)$ for the typed tree $b \in \skrix\{ k \}$
obtained by providing children for each vertex in generation $k-1$ of
$\widetilde{X} (v)$ independently according to the transition mechanism $\Q$.

With $\Q$ irreducible, it is easy to check that
any $a \in \skrix\{k-1\}$ can be reached by finitely many steps of
the transition mechanism $\widetilde{\Q}$ for $\widetilde{X}$
starting at any $b \in \skrix_{k-1}\{k-1\}$.
Further, if $b \in \skrix_l \{k-1\}$
for some $l < k-1$, then $\widetilde{\Q}\{\cdot\,|\,b\}$
is supported by $\bigcup_{n=0}^{N_0} \{n\} \times \skrix \{l-1\}^n$,
implying that $\widetilde{A}(a,b)=0$ whenever $a \in \skrix_m \{k-1\}$
for some $m \geq l$.
Consequently, $\widetilde{\Q}$ is weakly irreducible on $\skrix\{k-1\}$.
Let $\mu_0$ denote the Perron-Frobenius eigenvector
of the irreducible matrix $A$, normalized to be a strictly positive
probability vector on $\skrix$. Then,
$\mu_l = \mu_{l-1} \otimes_1 \Q$ for $l \geq 1$ are strictly positive
probability vectors on $\skrix \{l\}$, such that
$\mu_l \circ \pi^{-1}_{l,l-1} = \mu_{l-1}$ for all $l \geq 1$. Moreover,
with $\mu_0$ the right eigenvector corresponding to the eigenvalue $1$
of the matrix $A$, it follows by induction on $l \geq 1$ that $\mu_l$ are
shift-invariant on $\skrix(l)$. In particular,
for any $a \in \skrix\{k-1\}$,
$$\sum_{b \in \skrix\{k-1\}} \widetilde{A}(a,b) \mu_{k-1}(b) :=
\sum_{\heap{b\in\skrix\{k-1\}}{c\in\skrix\{k-1\}^*}}
\widetilde{m}(a,c) \widetilde{\Q}(c\,|\,b) \mu_{k-1}(b)
= \sum_{\bar{b} \in \skrix\{k\} } m_k(a,\bar{b})\,
\mu_{k-1} \otimes_1 \Q (\bar{b}) = \mu_{k-1}(a) \,.$$
With $\mu_{k-1}$ a strictly positive
right eigenvector for the eigenvalue $1$ and the matrix $\widetilde{A}$, we
see that $\widetilde{\Q}$ is also critical. Consequently,
we have from Theorem~\ref{general} that $M_{\widetilde{X}}$
satisfy the large deviation principle in
$\skrim(\skrix\{k-1\} \times \skrix\{k-1\}^*)$ with the good
rate function $\widetilde{J} (\cdot)$ corresponding to $\widetilde{\Q}$.
For each $\nu_1 \in \skrim(\skrix\{k-1\})$ the measure
$\nu_1 \circ \widetilde{\Q}$ is supported on
the closed (finite) set $\skrii(\skrix\{k\})$.  Consequently,
$M_{\widetilde{X}}$ is supported on $\skrii(\skrix\{k\})$ as is any
$\nu$ for which $\widetilde{J} (\nu) <\infty$, allowing us to
restrict this large deviation principle to $\skrim(\skrii(\skrix\{k\}))$.
Identifying $\skrim(\skrii(\skrix\{k\}))$ with
$\skrim(\skrix\{k\})$ via the mapping
$\mu=\nu \circ \skrii$,
the law of $M_{\widetilde{X}}$ is exactly mapped to that of $M_X^k$.
Moreover, $\nu \in \skrim(\skrii(\skrix\{k\}))$ is shift-invariant
if and only if $\mu$ is shift-invariant on $\skrix(k)$ as defined
in \eqref{disp0}, with $\nu_1 = \mu \circ \pi_{k,k-1}^{-1}$
and $(\nu_1 \otimes \widetilde{\Q}) \circ \skrii
 =(\mu \circ \pi_{k,k-1}^{-1}) \otimes_1 \Q$. This leads to the
large deviation principle for $M^k_X$ with the
good rate function $J_k(\cdot)$, restricted to $\skrim(\skrix\{k\})$.

To complete the proof it suffices to check that any shift-invariant
measure $\mu\in\skrim(\skrix(k))$ with $\mu \ll \mu \circ \pi_{k,k-1}^{-1}
\otimes_1 \Q$ in $\skrim(\skrix(k))$ is supported by $\skrix \{k\}$.
To this end, fix a shift-invariant $\mu$ in $\skrim(\skrix(k))$
and note that $\int N[m] d\mu =1$ for $m=1,\ldots,k$. Hence
we can associate shifted probability measures
$S^m (\mu) \in \skrim(\skrix(k-m))$ with $\mu$ such that $S^0(\mu)=\mu$,
$S^m (\mu) = S (S^{m-1}(\mu))$ for $m=1,\ldots,k$, and
$S(\mu)$ is defined 
as in \eqref{s-def}. The shift-invariance of $\mu$ implies that
$S^m (\mu) \circ \pi^{-1}_{k-m,1}$ is independent of $m=0,\ldots,k-1$.
Recall that the measure
$S^{k-1}(\mu)$ of each $(a,c) \in \skrix(1)$ is the expectation
under $\mu$ of the number of vertices of generation $k-1$ of the tree
whose type is $a \in \skrix$ and which have offspring $c \in \skrix^*$.
Our assumption that $\mu \ll \mu \circ \pi_{k,k-1}^{-1} \otimes_1 \Q$
thus implies that the support of $S^{k-1}(\mu)$ is a subset
of the support of $\mu_1$, which is $\skrix\{1\}$.
Consequently, $S^m (\mu) \circ \pi^{-1}_{k-m,1}$
are supported by $\skrix \{ 1 \}$ for all $m=0,\ldots,k-1$,
which implies that $\mu$ is supported by $\skrix \{k\}$ as claimed.
\end{proof}

To move from the empirical $k$-generation offspring measures $M_X^k$
to the empirical subtree measure $T_X$ we use the Dawson-G\"artner
theorem, see e.g.~\cite[Theorem 4.6.1]{DZ98}. Note that the spaces
$\skrix(k)$ and the canonical projections $\pi_{k,l}$, $k\ge l$,
form a projective system
of Polish spaces
and that the projective limit coincides with
the Polish space $\bar{\skrix}$.

Similarly, the probability measures on $\skrix(k)$
with the projections $\pi^*_{k,l}$ defined by $\pi^*_{k,l}(\mu)=
\mu\circ\pi_{k,l}^{-1}$ form a projective system and the projective limit
is the
Polish space $\skrim(\bar{\skrix})$
described before Theorem~\ref{process} and the canonical
projections $\pi^*_{k}
: \skrim(\bar{\skrix}) \to \skrim(\skrix(k))
$ can be defined by $\pi^*_{k}(\mu)=
\mu\circ\pi_{k}^{-1}$. Details follow from an argument similar to
the one given in \cite[Lemma 6.5.14]{DZ98}. Recalling that
$M_X^k=T_X\circ \pi_k^{-1}$, the Dawson-G\"artner theorem yields the
following corollary of Lemma \ref{kstep}
(see for example \cite[Corollary 6.5.15]{DZ98} for a similar derivation).

\begin{cor}\label{process-ldp}
Suppose that $X$ is an irreducible, critical multitype Galton-Watson
tree
with uniformly bounded offspring sizes,
conditioned to have exactly $n$ vertices. Then, for $n\to\infty$,
the empirical subtree measure~$T_X$ satisfies a large deviation
principle in $\skrim(\bar{\skrix})$ with speed $n$ and
convex,
good rate function
$$\tilde{K}(\mu)=\sup_{k\ge 1} J_k\big( \mu\circ \pi_k^{-1}\big).$$
\end{cor}

To complete the proof of Theorem \ref{process}
it just remains to show that $\tilde{K}(\cdot)=K(\cdot)$. For this
purpose first assume that $\mu
\in \skrim(\bar{\skrix})$
is shift-invariant. Then, for each $k\ge 1$
and $a\in \skrix(k-1)$,
\begin{align*}
\big( \mu\circ\pi_{k}^{-1} \big) \circ\pi_{k,k-1}^{-1}(a) = &
\, S(\mu)\circ\pi_{k-1}^{-1}(a)
=   \int d\mu(X) \,\sum_{i=1}^N \delta_{X^{v_i}_{k-1}}(a) \\
= & \int d\mu(X) \, m_k(a, \pi_kX)  =  \sum_{b\in\skrix(k)}
\mu\circ\pi_k^{-1}(b) \, m_k(a,b).
\end{align*}
In other words, for each $k\ge 1$, the measure $\mu\circ\pi_{k}^{-1}$ is
shift-invariant
in $\skrim(\skrix(k))$.
Conversely, if $\mu\circ \pi^{-1}_k$ is shift-invariant
in $\skrim(\skrix(k))$
for
every $k\ge 1$, the same calculation shows that $\mu=S(\mu)$ on the collection
of sets of the form $\pi^{-1}_k(A)$ for any $k\ge 1$ and $A\subset\skrix(k)$.
As this collection of sets is closed
under finite intersections and it generates the Borel $\sigma$-field
on $\bar{\skrix}$, we infer that $\mu$ itself
is shift-invariant.

Recall the definition of the projections
$\mathfrak{p}_0, \mathfrak{p}_1$ for backward trees.
For the proof of Theorem~\ref{process} it only remains to verify
the following lemma.

\begin{lemma}\label{lem-Kiden}
For every shift-invariant probability measure $\mu$ on $\bar{\skrix}$
we have
\begin{equation}\label{identi}
H\big( \mu^*\circ \mathfrak{p}_1^{-1} \, \big\| \,
\mu^*\circ \mathfrak{p}_{0}^{-1}\otimes   \Q \big)
= \sup_{k\ge 2} H\big( \mu\circ \pi_k^{-1} \, \big\| \,
\mu\circ \pi_{k-1}^{-1}\otimes_1 \Q \big).
\end{equation}
\end{lemma}

\begin{proof}
Define projections $\pi_k^j:\bar{\skrix}\to\skrix(k)$ as follows:
Order the vertices $v_1,v_2,\ldots$ in generation $k-1$ of $x\in\bar{\skrix}$
from left to right, with $v_1$ the leftmost. The tree $\pi_k^j(x)$ is
obtained by removing all vertices in
generations exceeding $k$ and all vertices in generation $k$ whose
parent is some $v_l$, $l \geq j$. In particular, $\pi_k^1(x)=\pi_{k-1}(x)$
and $\pi_k^j(x)=\pi_k(x)$ for all $j>N[k-1](x)$. Let
$\mu\circ (\pi_k^j)^{-1}\otimes_j\Q$ denote the measure
obtained by sampling $X$ according to $\mu$ and
then independently adding offspring according to $\Q$
to each of the vertices $v_l$ for $l\geq j$ in generation $k-1$ of
$\pi_k^j(X)$. Observe that we
define this measure for \emph{all} $j$ and that in many cases no
vertices in generation $k$ are removed or added.
Assume first
that $\mu\circ \pi_k^{-1} \ll \mu\circ \pi_{k-1}^{-1}\otimes_1 \Q$.
Then, in case $\mu\circ \pi_k^{-1}(x)>0$ and $N[k-1](x) = n \geq 1$
we find that
\begin{equation}\label{disp1}
\frac{\mu\circ \pi_k^{-1}(x)}{\mu\circ \pi_{k-1}^{-1}\otimes_1 \Q(x)} =
\prod_{j=1}^{n} \frac{\mu\circ (\pi_k^{j+1})^{-1}\otimes_{j+1}\Q(x)}
{\mu\circ (\pi_k^j)^{-1}\otimes_{j}\Q(x)},
\end{equation}
with all the terms on the right hand side positive.
Recall the definition of the measure
$\mu_{k-1}=\mu^*\circ p_{k-1}^{-1}$ and
the projections $\mathfrak{p}_{0,k-1}$,
$\mathfrak{p}_{1,k-1}$ on $\skrix[k-1]$ 
and also recall that $(y,v)\in\skrix[k-1]$ denotes
the tree $y \in \bar{\skrix}$ with centre $v$ in
generation $k-1$ of $y$. Hence, for $1 \leq j \leq n$,
\begin{equation}\label{disp2}
\frac{\mu\circ (\pi_k^{j+1})^{-1}\otimes_{j+1}\Q(x)}
{\mu\circ (\pi_k^j)^{-1}\otimes_{j}\Q(x)}
=\frac{\mu_{k-1}\circ \mathfrak{p}_{1,k-1}^{-1}(\pi_k^j(x),v_j)}
{\mu_{k-1}\circ\mathfrak{p}_{0,k-1}^{-1}\otimes\Q(\pi_k^j(x),v_j)}
\end{equation}
with all terms positive. Note that if $N[k-1](x)=0$ then
$\mu\circ \pi_k^{-1}(x)=\mu\circ \pi_{k-1}^{-1}\otimes_1 \Q(x)$,
whereas if $y=\pi_k^j(x)$ with $N[k-1](x)=n >0$ then
$N[k-1](y)=n$ and
$\mu\circ (\pi_k^j)^{-1} (y) =\mu\circ \mathfrak{p}_{1,k-1}^{-1}(y,v_j)$
for any $1 \leq j \leq n$.
Hence, \eqref{disp1} and \eqref{disp2} imply that
\begin{align}
H\big( \mu\circ \pi_k^{-1} \, \big\| \,
\mu\circ \pi_{k-1}^{-1}\otimes_1 \Q \big) &=
\sum_{x\in\skrix(k)} \sum_{j=1}^{N[k-1](x)} \,
\mu\circ\pi_k^{-1}(x)
\log \Big(
\frac{\mu_{k-1}\circ \mathfrak{p}_{1,k-1}^{-1}(\pi_k^j(x),v_j)}
{\mu_{k-1}\circ\mathfrak{p}_{0,k-1}^{-1}\otimes\Q(\pi_k^j(x),v_j)}
\Big) \nonumber \\
&= H\big( \mu_{k-1} \circ \mathfrak{p}_{1,k-1}^{-1} \, \big\| \,
\mu_{k-1} \circ \mathfrak{p}_{0,k-1}^{-1} \otimes \Q \big) .
\label{disp3}
\end{align}
Finally, note that 
$$\mu \circ \pi_k^{-1} (x) >0 \mbox{ and }
\mu\circ \pi_{k-1}^{-1}\otimes_1 \Q (x)=0 \mbox{  for some
$x \in \skrix(k)$},$$
if and only if there exists $1 \leq j \leq N[k-1](x)$ such that
$$\mu_{k-1} \circ \mathfrak{p}_{1,k-1}^{-1} (\pi_k^j(x),v_j) >0
\mbox{ and }
\mu_{k-1} \circ \mathfrak{p}_{0,k-1}^{-1} \otimes \Q (\pi_k^j(x),v_j) =0.$$
Consequently, $\mu\circ \pi_k^{-1} \ll \mu\circ \pi_{k-1}^{-1}\otimes_1 \Q$
if and only if $\mu_{k-1} \circ \mathfrak{p}_{1,k-1}^{-1} \ll
\mu_{k-1} \circ \mathfrak{p}_{0,k-1}^{-1} \otimes \Q$, with
\eqref{disp3} holding for any shift-invariant $\mu \in \skrim(\bar{\skrix})$
and $k \geq 2$. By the identities
$p_k\circ \mathfrak{p}_0=\mathfrak{p}_{0,k}\circ p_k$
and $p_k \circ \mathfrak{p}_1=\mathfrak{p}_{1,k}\circ p_k$
this amounts to
\begin{equation}
H\big( \mu\circ \pi_k^{-1} \, \big\| \, \mu\circ \pi_{k-1}^{-1}\otimes_1
\Q \big) = H\big(
\mu^*\circ \mathfrak{p}_1^{-1}\circ p_{k-1}^{-1} \, \big\| \, (\mu^*\circ
\mathfrak{p}_{0}^{-1}\otimes\Q)\circ p_{k-1}^{-1} \big) \,.
\label{key}
\end{equation}
The variational characterization of the relative entropy states that, for two
probability measures $\nu_1, \nu_2$ on
the Polish space
$\underline{\skrix}$,
$$H\big( \nu_1\circ p_{k}^{-1} \, \big\| \,\nu_2\circ p_{k}^{-1} \big)
=\sup_{\phi\in C_b(\skrix[k])} \bigg\{ \int_{\underline{\skrix}}
\phi\circ p_{k} \, d\nu_1
- \log \int_{\underline{\skrix}} e^{\phi\circ p_{k}} \, d\nu_2\bigg\},$$
where $C_b(\skrix[k])$ is the set of continuous, bounded
functions on $\skrix[k]$ (see for example \cite[Lemma 6.2.13]{DZ98}).
Obviously, this expression is increasing in $k$ and by the same
representation it is bounded by
$H(\nu_1\,\|\, \nu_2)$, which together with (\ref{key}) shows
that the left hand side of \eqref{identi} is at least as large as its
right hand side.

Conversely, for any continuous bounded function
$\phi:\underline{\skrix}\to\R$
and $\eps>0$ there exists a uniformly continuous function
$\psi:\underline{\skrix}\to\R$ such that
$$\Big| \log \int_{\underline{\skrix}} e^{\phi} \, d\nu_2
-\log \int_{\underline{\skrix}} e^{\psi} \, d\nu_2\Big|
<\eps \quad \mbox{ and } \quad
\Big| \int_{\underline{\skrix}} \phi \, d\nu_1
-\int_{\underline{\skrix}} \psi \, d\nu_1\Big| <\eps.$$
Moreover,
with $\underline{\skrix}$ being the projective limit of $\skrix[k]$,
we can find a $k\ge 1$ and a
continuous, bounded
function $\psi^k:\skrix[k]\to\R$ such that
$|\psi^k\circ p_{k}(x)-\psi(x)|<\eps$ for all $x\in\underline{\skrix}$.
Hence $$\sup _{k\ge 2} \sup_{\phi\in C_b(\skrix[k])}
\bigg\{ \int_{\underline{\skrix}} \phi\circ p_{k} \, d\nu_1
- \log \int_{\underline{\skrix}} e^{\phi\circ p_{k}} \, d\nu_2\bigg\}
\ge \sup_{\phi\in C_b(\underline{\skrix})}
\bigg\{ \int_{\underline{\skrix}} \phi \, d\nu_1
- \log \int_{\underline{\skrix}} e^{\phi} \, d\nu_2\bigg\},$$
which together with (\ref{key}) shows
that the right hand side of \eqref{identi} is at least as large as its
left hand side. This completes the proof of the lemma.
\end{proof}


\bigskip


{\it AMIR DEMBO, Department of Mathematics, Stanford University\\
Stanford, CA 94305, USA.

PETER M\"ORTERS, Department of Mathematical Sciences, University of Bath\\
Bath BA2 7AY, United Kingdom.

SCOTT SHEFFIELD, Microsoft Research\\
One Microsoft Way, Redmond WA 98052, USA.}

\end{document}